# Bernstein polynomials: a bibliometric data analysis since the year 1949 based on the Scopus database

**Rushan Ziatdinov** [1,*]

[1] Department of Industrial Engineering, College of Engineering, Keimyung University, 704-701 Daegu, Republic of Korea
* Correspondence: ziatdinov@kmu.ac.kr, ziatdinov.rushan@gmail.com, URL: https://www.ziatdinov-lab.com/

**Abstract:** It's hard to imagine human life in the digital and AI age without polynomials because they are everywhere but mostly invisible to ordinary people: in data trends, on computer screens, in the shapes around us, and in the very fabric of technology. One of these, the simple but elegant Bernstein polynomials, was discovered by a scientist from the Russian Empire, Sergei Bernstein, in 1912 and plays a central role in mathematical analysis, computational and applied mathematics, geometric modelling, computer-aided geometric design, computer graphics and other areas of science and engineering. They have been the subject of much research for over a hundred years. However, no work has carried out database-derived research analysis, such as bibliometric, keyword or network analysis, or more generally, data analysis of manuscript data related to Bernstein polynomials extracted from digital academic databases. This work, which appears to be the first-ever attempt at the bibliometric data analysis of Bernstein polynomials, aims to fill this gap and open researchers' eyes to potentially new or underexplored areas of mathematics and engineering where Bernstein polynomials may one day be used to make discoveries. The results may be helpful to academics researching Bernstein polynomials and looking for potential applications, collaborators, supervisors, funding or journals to publish in.

---

*"Today all mathematicians and physicists agree that the field of applications for mathematics knows no limits except those of knowledge itself."*
*Sergei Bernstein*

## 1. Introduction

It's hard to imagine human life in the digital and AI age without polynomials, as they are everywhere around us but often invisible to most people: in data trends, on computer screens, in the shapes around us, and the very fabric of technology. Polynomials drive computer animation in movies, optimise routes in navigation systems such as GPS, model complex physical phenomena from planetary motion to the spread of infectious diseases, and are part of models of many complex systems that help us better understand the world around us. They're at the heart of the algorithms that process images, the central part of the geometric kernels in CAD systems. They're used to predict weather conditions and compress the music we listen to daily. In engineering, they define structural designs and stress analysis; in economics, they model growth and fluctuations; in machine learning, they help build predictive models. Polynomials are not just equations made up of mathematical symbols - they are the silent architects of the modern digital world, shaping how we interact with technology, interpret data, and design innovation contours.

One of these polynomials is the simple but elegant Bernstein polynomial, which has been a pivotal concept in mathematics since it was introduced by Sergei Bernstein in 1912. These polynomials have the following form

$$B_n(x) = \sum_{\nu=0}^{n} \beta_\nu b_{\nu,n}(x)$$

where $\beta_\nu$ is a Bernstein coefficient, and $(n+1)$ Bernstein basis polynomials of degree $n$ are defined by

$$b_{\nu,n}(x) = \binom{n}{\nu} x^\nu (1-x)^{n-\nu}, \nu = 0, \dots, n,$$

where $\binom{n}{\nu}$ is a binomial coefficient.

Initially used in a constructive proof for the Weierstrass approximation theorem, these simple but elegant polynomials have evolved into a fundamental tool for approximation theory, numerical analysis, and computer-aided geometric design (CAGD) [1]. Their simple yet profound nature has allowed them to remain an active research topic in CAGD for more than one hundred years, stimulating further exploration and various applications in engineering.

The enduring significance of Bernstein polynomials stems from their simplicity, elegance, and versatility. Their key properties, such as non-negativity, partition of unity, and shape-preserving characteristics, make them indispensable for mathematical modelling, curve fitting, and optimisation problems [2]. Their association with Bézier curves and surfaces in the form of basis functions has also made them a cornerstone in computer graphics, computer-aided design (CAD), and VR, enabling the creation of visually appealing and mathematically precise geometric shapes vital for digital manufacturing and 3D printing.

Sergei Bernstein died in 1968 at the age of 88, which means that polynomials with his name were used (at least in international journals) for about twenty years of his last years. Out of a total of 16,392 Scopus-indexed documents containing the keywords "Bézier surface" or "Bézier surface", "Bernstein-Bézier curve" or "Bernstein-Bézier surface" are only mentioned in 416 documents in all search fields, i.e. only 2.57%. From one point of view, this is a pity, but from another point of view, these days, who mentions the software or programming language used for their discoveries in the title of a manuscript or the name of a new theorem or method? The various tools we use in research have often become the silent architects of scientific progress.

The applications of Bernstein polynomials and their generalisations have expanded across diverse disciplines, including recognition of human speech [3], probability [4], and the theory of desirable gambles [5]. This growth has been accompanied by extensive research, continuously enhancing our understanding of these polynomials and their properties and reflecting their ongoing relevance. Over the past century, the research on Bernstein polynomials has experienced significant growth and development [6].

This manuscript presents a pioneering bibliometric data analysis of Bernstein polynomials based on RIS data file download from the Scopus database (http://scopus.com/), analysing publications since 1949. The aim is to examine the evolution of research topics, identify key contributors and organisations, funding agencies, and influential works, and search for areas where Bernstein polynomials are rarely mentioned. This analysis aims to provide a valuable resource that encapsulates Bernstein polynomials' historical progress and future potential. The Scopus search query TITLE-ABS-KEY(polynomial AND "bibliometric analysis") returns 23 papers, but "polynomial" is not in the title of these papers, which means that our work may be the first-ever attempt at bibliometric analysis of polynomials, at least based on papers indexed in the Scopus database.[1] of polynomials, at least based on papers indexed in the Scopus database.

## 2. Materials and Methods

### *2.1. Data Processing*

A search was conducted for documents (excluding preprints, patents, and secondary documents) with the keyword "Bernstein polynomial" in all fields. Scopus, one of the largest web bibliometric databases, was used to obtain the dataset. Several keywords were used to find weak connections to Bernstein polynomials and to attract academics for research on Bernstein polynomials applied to found areas. Our study excludes documents that may use Bernstein polynomials under other names, as they were called before the term "Bernstein polynomials" became widely accepted in the academic community.

### *2.2. Analysis Procedure and Software*

To carry out a bibliometric analysis, we followed the generally accepted steps: (1) study or research design, (2) data collection, (3) data analysis, (4) data visualisation, and (5) interpretation of results. We used VOSviewer v. 1.6.20 software and various AI chatbots such as ChatGPT, Llama, Mixtral, and Perplexity to analyse and visualise our data. The simple statistics, including the most productive authors, organisations, and others, were extracted from the Scopus database, and the relative percentages were calculated in Microsoft Excel. To carry out the co-citation analysis, we used the VOSviewer[2] software, which is used in bibliometric analysis to map networks and evaluate the strength of links.

We used the standard 32-bit colour maps available in VOSviewer to visualise the networks. The colour maps used facilitate the presentation of large amounts of bibliometric data and allow clear visualisation of network structures. However, users can easily customise and extend this colour mapping functionality by designing their colour maps for more accurate colour representation. This flexibility allows the creation of more nuanced and detailed visualisations, particularly for huge network datasets.

---

[1] Authors should not confuse ***bibliometric analysis*** with ***bibliographic analysis***, which are related but distinct concepts in the study of academic literature and research output.

[2] https://www.vosviewer.com/

When exporting visualisations in VOSviewer, the software adds its logo in the lower left corner, but it's missing on screenshots. Again, we want to emphasise that all visualisations have been made with the software mentioned. As no copyright information exists on their website, adding the logo could be seen as an advertising method.

## 3. Basic Statistics from Scopus Database

Basic statistics, such as subject area, source title, etc., can be copied or scraped directly from the Scopus search results page. As of 23 January 2025, the keyword "Bernstein polynomial" is mentioned in 2371 documents' titles, abstracts or keywords, and the search within all fields returns 9310 documents starting from 1949. Such a significant difference may indicate that many papers cite "Bernstein polynomials" without directly researching them but using them in various fields. As we are looking for multiple applications of these polynomials, we would consider the dataset with 9310 documents, including 77 articles in press, for our research. Of these, 84.9% are articles, 11% conference papers, 1.8% book chapters, 1.1% reviews, 0.8% books, 0.2% notes, 0.1% editorials and 0.3% other sources. Interestingly, only two retracted papers among the documents show a high academic culture level in this field.

Most of the fields where the Bernstein polynomial is mentioned are shown in Figure 1, where the most popular fields are Mathematics, Computer Science, Engineering, and Physics & Astronomy. Fields such as Materials Science, Chemistry, Earth and Planetary Sciences, Energy, Biochemistry, Genetics and Molecular Biology, Agricultural and Biological Sciences, Multidisciplinary, Environmental Science, Economics, Econometrics and Finance, Social Sciences, Chemical Engineering, Medicine, Business, Management and Accounting, Neuroscience, Immunology and Microbiology, Arts and Humanities, Health Professions, Pharmacology, Toxicology and Pharmaceutics, Psychology, Nursing, Dentistry have less than 3% of all documents for each field.

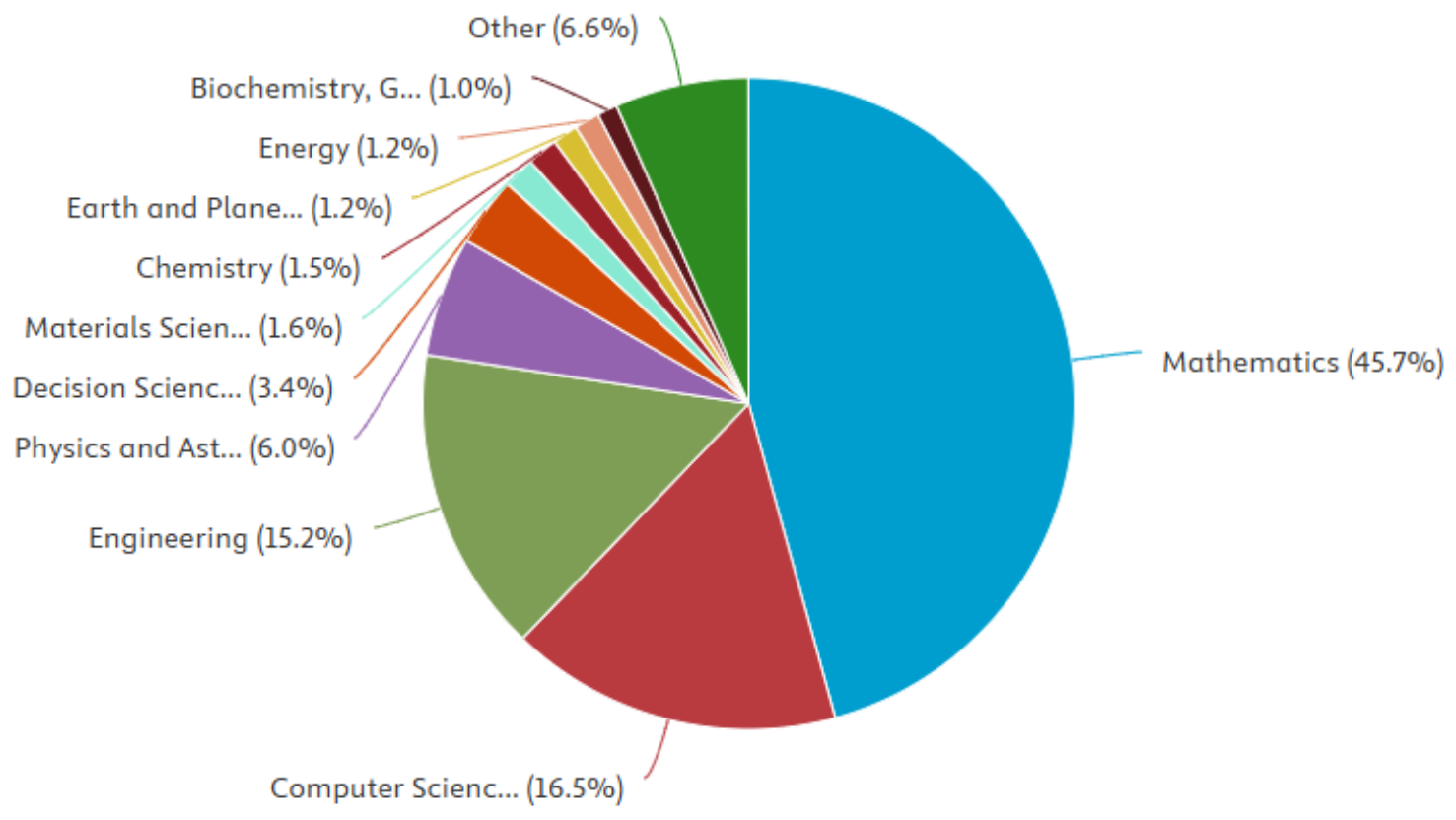

**Figure 1.** A pie chart from Scopus shows the percentage of total mentions of the keyword "Bernstein polynomial" in different subject areas.

Figure 2 shows the dynamics of annual publications and the twenty most popular keywords since 1949. From around 2005, we can observe a significant growth in the number of documents and, consequently, in the number of keywords.

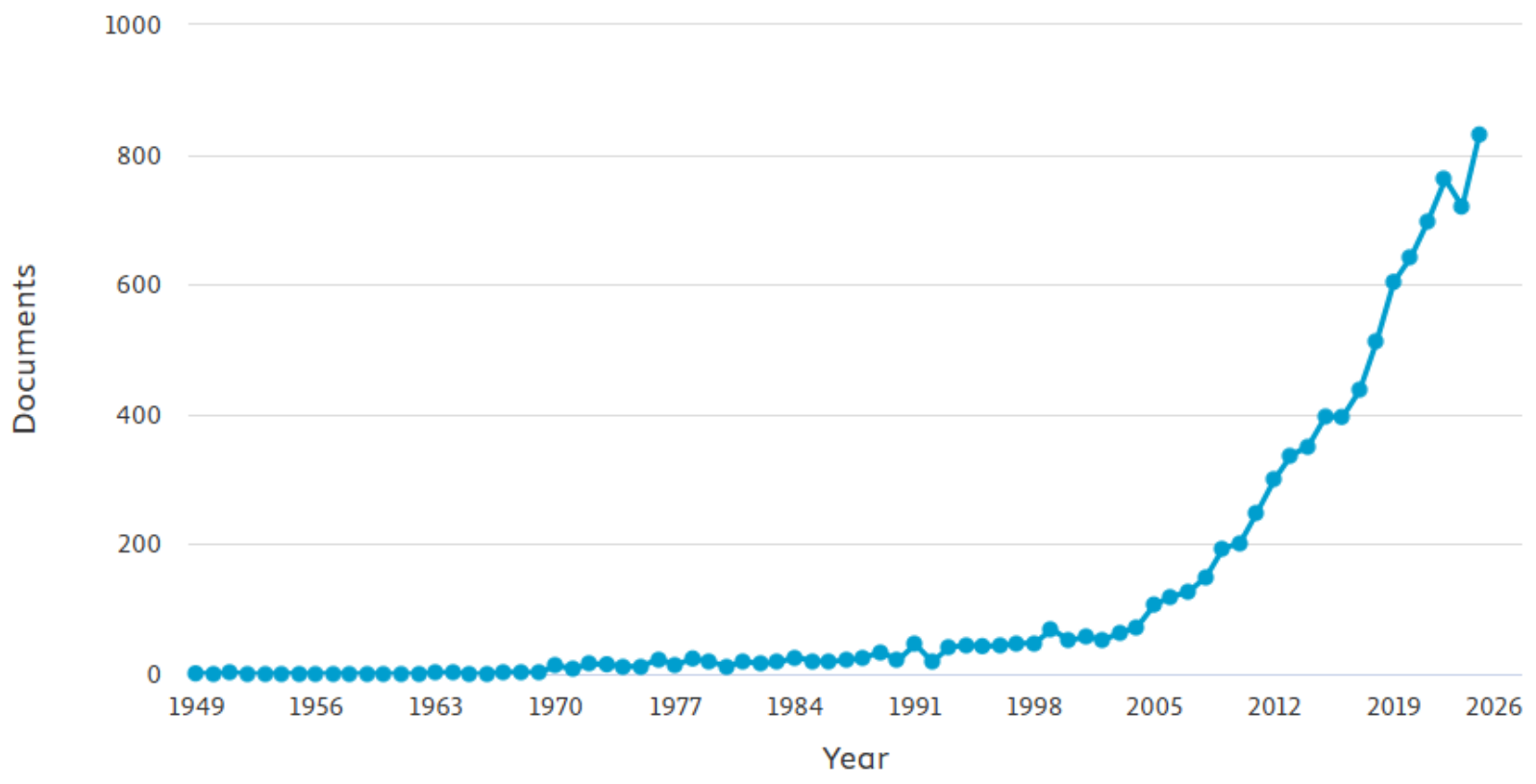

**(a)**

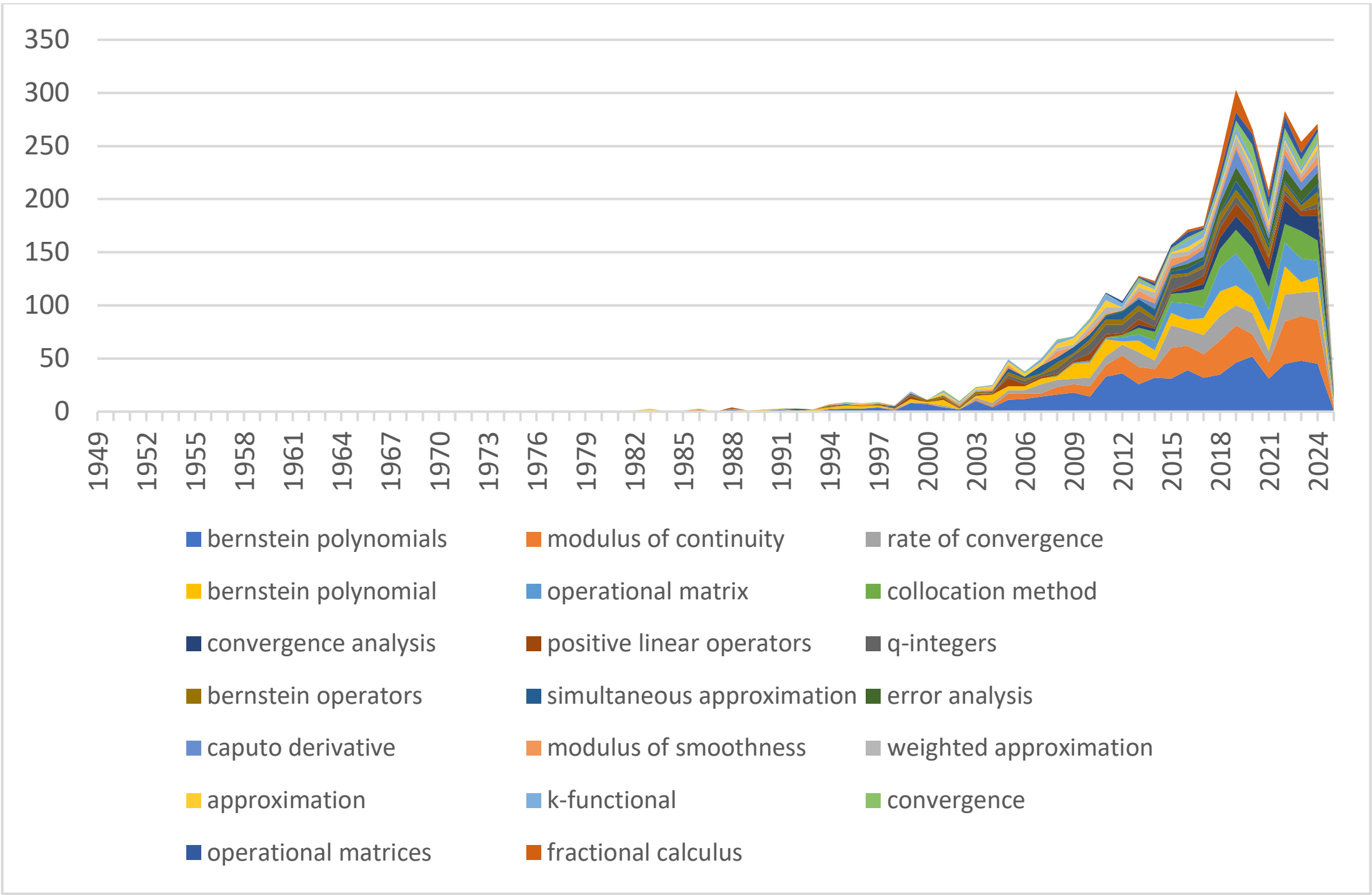

**(b)**

**Figure 2.** (a) Number of documents since 1949. The year 2025 is not shown here, and 70 publications have been indexed in Scopus. (b) The evolution of the twenty most popular keywords in manuscripts. The evolution of keywords with ten or more occurrences can be found online at https://youtu.be/DQBZHmngCNo

The 100 most popular sources, mainly containing prestigious journals, account for 4060 or 43.6% of the 9310 documents in the database (see Table 1). Most of these are mathematics, computer, and applied sciences-related periodicals.

**Table 1.** The list of the most popular source titles where the keyword "Bernstein polynomial" appears in all search fields.

| Source name | Number of documents |
|---|---|
| Journal of Approximation Theory | 240 |
| Applied Mathematics and Computation | 205 |
| Journal of Computational and Applied Mathematics | 165 |
| Mathematical Methods in the Applied Sciences | 138 |
| Journal of Inequalities and Applications | 112 |
| Journal of Mathematical Analysis and Applications | 105 |
| Filomat | 101 |
| Computer Aided Geometric Design | 97 |
| Mathematics (MDPI) | 95 |
| Journal of Computational Analysis and Applications | 94 |
| Advances in Difference Equations | 91 |
| Results in Mathematics | 89 |
| Lecture Notes in Computer Science, including Subseries Lecture Notes in Artificial Intelligence and Lecture Notes in Bioinformatics | 78 |
| AIP Conference Proceedings | 76 |
| Symmetry | 66 |
| Abstract and Applied Analysis | 61 |
| Computational and Applied Mathematics | 61 |
| International Journal of Applied and Computational Mathematics | 55 |
| Computers and Mathematics with Applications | 54 |
| Numerical Algorithms | 54 |
| Advanced Studies in Contemporary Mathematics Kyungshang | 51 |
| Numerical Functional Analysis and Optimization | 46 |
| Mediterranean Journal of Mathematics | 44 |
| Applied Numerical Mathematics | 42 |
| Chaos Solitons and Fractals | 40 |
| Computational Statistics and Data Analysis | 40 |
| Proceedings of The Jangjeon Mathematical Society | 40 |
| Mathematics and Computers in Simulation | 39 |
| Demonstratio Mathematica | 36 |
| Aims Mathematics | 35 |
| Applied Mathematical Sciences | 35 |
| Acta Mathematica Hungarica | 34 |
| Constructive Approximation | 33 |
| Iranian Journal of Science and Technology Transaction A Science | 33 |
| Carpathian Journal of Mathematics | 32 |
| Springer Proceedings in Mathematics and Statistics | 32 |
| Fractal and Fractional | 31 |
| Engineering with Computers | 30 |
| Advances in Computational Mathematics | 28 |
| Alexandria Engineering Journal | 28 |
| Applied Mathematical Modelling | 28 |
| IEEE Transactions on Antennas and Propagation | 28 |
| Applied Mathematics Letters | 27 |
| Communications in Statistics Theory and Methods | 27 |
| IEEE Access | 27 |
| International Journal of Mathematical Analysis | 27 |
| Journal of Nonparametric Statistics | 27 |
| Turkish Journal of Mathematics | 27 |
| Applied Mathematics and Information Sciences | 26 |
| Approximation Theory and Its Applications | 26 |
| Journal of Physics Conference Series | 26 |
| Mathematical Sciences | 26 |
| Revista De La Real Academia De Ciencias Exactas Fisicas Y Naturales Serie A Matematicas | 26 |
| Journal of Function Spaces | 25 |
| Mathematical Foundations of Computing | 25 |
| International Journal of Computer Mathematics | 24 |
| Journal of Vibration and Control | 24 |
| Journal of Statistical Planning and Inference | 24 |
| Springer Optimization and Its Applications | 24 |
| Computational Methods for Differential Equations | 23 |
| Computer Methods in Applied Mechanics and Engineering | 23 |
| Journal of Mathematical Inequalities | 23 |
| Soft Computing | 23 |
| Calcolo | 22 |
| Fuzzy Sets and Systems | 22 |
| IEEE Transactions on Signal Processing | 22 |
| Journal of Multivariate Analysis | 22 |
| Positivity | 22 |
| International Journal of Pure and Applied Mathematics | 21 |
| Mathematical and Computer Modelling | 21 |
| Mathematical Problems in Engineering | 21 |
| Rocky Mountain Journal of Mathematics | 21 |
| Statistics and Probability Letters | 21 |
| Studia Scientiarum Mathematicarum Hungarica | 21 |
| Applied Mathematics | 20 |
| CMES Computer Modeling in Engineering and Sciences | 20 |
| Electronic Journal of Statistics | 20 |
| International Journal of Mathematics and Mathematical Sciences | 20 |
| Numerical Methods for Partial Differential Equations | 20 |
| Hacettepe Journal of Mathematics and Statistics | 19 |
| IEEE Signal Processing Letters | 19 |
| Journal of Computational Physics | 19 |
| Journal of The American Statistical Association | 19 |
| Advanced Studies in Theoretical Physics | 18 |
| Communications in Statistics Simulation and Computation | 18 |
| Georgian Mathematical Journal | 18 |
| Journal of Numerical Analysis and Approximation Theory | 18 |
| Mathematica Slovaca | 18 |
| Miskolc Mathematical Notes | 18 |
| Nonlinear Dynamics | 18 |
| Proceedings IEEE International Symposium on Circuits and Systems | 18 |
| Proceedings of The American Mathematical Society | 18 |
| Studia Universitatis Babes Bolyai Mathematica | 18 |
| Annals of Statistics | 17 |
| Discrete Dynamics in Nature and Society | 17 |
| Journal of Applied Mathematics and Computing | 17 |
| Journal of Intelligent and Fuzzy Systems | 17 |
| Axioms | 16 |
| Energy | 16 |
| Fractals | 16 |

Figure 3 shows the dynamics of the number of mentions of the keyword "Bernstein polynomial" in the four most popular journals from Table 1. The number of mentions of the keyword in Applied Mathematics and Computation has significantly decreased since 2015, while the number of mentions of the keyword "Bernstein polynomial" in the journal Mathematical Methods in the Applied Sciences has increased. Such changes could be related to the increased difficulty of publishing in the mentioned journal, changes in editorial policy or the impact factor, preferences of academics, etc.

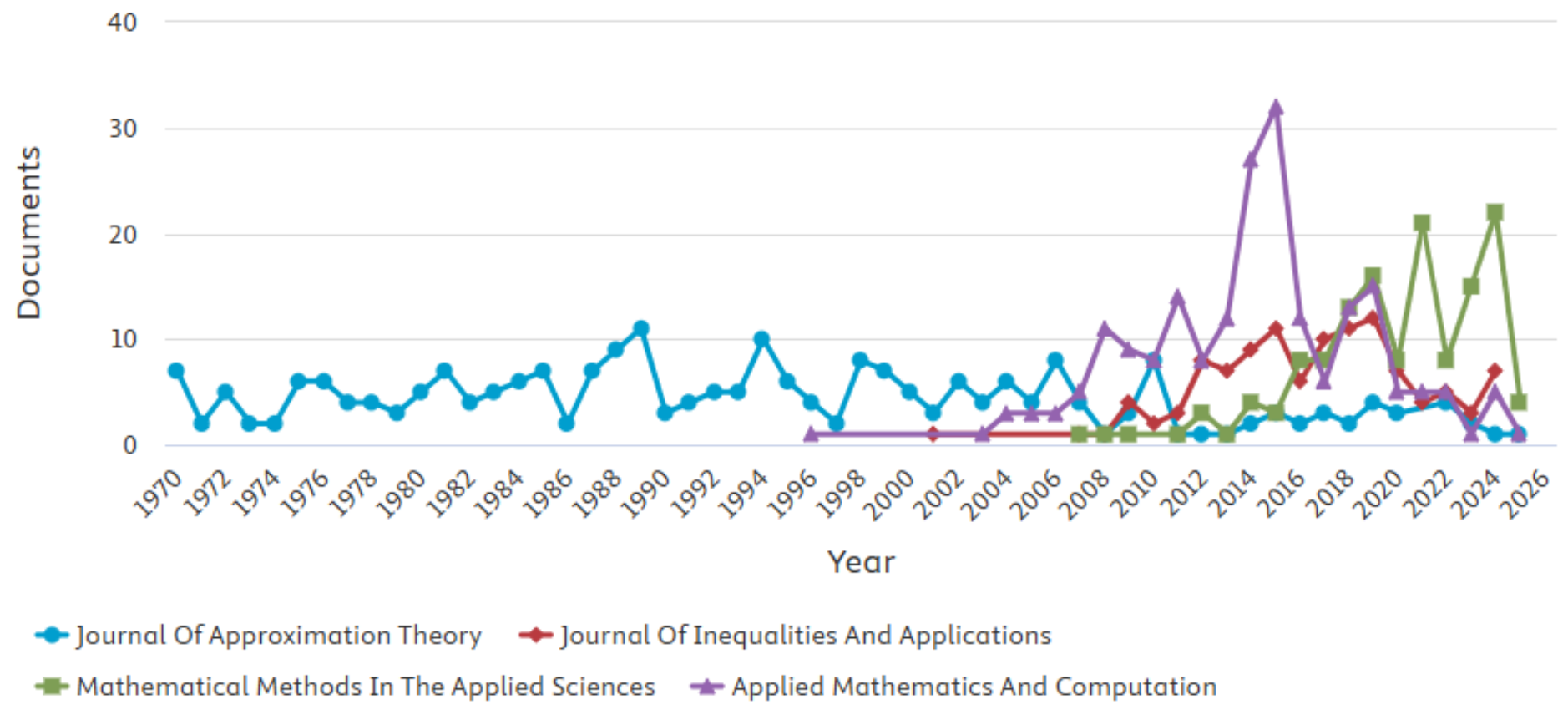

**Figure 3.** The dynamics of the number of mentions of the keyword "Bernstein polynomial" in the four most popular journals from Table 1.

Table 2 shows the list of authors' countries where documents with the keywords "Bernstein polynomial" were found in one of the search fields. As can be seen, China, USA and India dominate the list. Russia and Libya are mentioned twice in the Scopus database under different names for which the corresponding number of papers has been added. Some papers were not identified by geographical location, but their total share is only 1.3%. Russia's 19th place is surprising, as Bernstein polynomials were discovered there. However, this may be due to the authors' preference to publish their related research in Russian sources that are not indexed in academic databases. E-library (https://elibrary.ru/), a central Russian academic database, mentions only 67 documents for "полином Бернштейна" (singular form), 48 documents for "многочлен Бернштейна" (singular form), 101 documents with "многочлены Бернштейна" (plural form) and 182 documents with "полиномы Бернштейна" (plural form) in all fields, including the full text of a document.

**Table 2.** The list of authors' countries.

| Country/Territory | Documents | Percentage |
|---|---|---|
| China | 1530 | 16.43% |
| United States | 1371 | 14.73% |
| India | 1113 | 11.95% |
| Turkey | 935 | 10.04% |
| Iran | 899 | 9.66% |
| Romania | 509 | 5.47% |
| South Korea | 508 | 5.46% |
| Germany | 494 | 5.31% |
| Saudi Arabia | 389 | 4.18% |
| Italy | 373 | 4.01% |
| France | 299 | 3.21% |
| Canada | 278 | 2.99% |
| United Kingdom | 268 | 2.88% |
| Spain | 259 | 2.78% |
| Egypt | 195 | 2.09% |
| Taiwan | 188 | 2.02% |
| Pakistan | 165 | 1.77% |
| Japan | 161 | 1.73% |
| Russia | 151 | 1.62% |
| Poland | 143 | 1.54% |
| Malaysia | 135 | 1.45% |
| Iraq | 127 | 1.36% |
| Australia | 126 | 1.35% |
| Jordan | 92 | 0.99% |
| South Africa | 76 | 0.82% |
| Switzerland | 74 | 0.79% |
| Netherlands | 71 | 0.76% |
| Bulgaria | 70 | 0.75% |
| Thailand | 70 | 0.75% |
| Portugal | 69 | 0.74% |
| Israel | 67 | 0.72% |
| Hong Kong | 65 | 0.70% |
| Singapore | 65 | 0.70% |
| Belgium | 64 | 0.69% |
| Mexico | 64 | 0.69% |
| Algeria | 60 | 0.64% |
| Hungary | 54 | 0.58% |
| Brazil | 50 | 0.54% |
| Viet Nam | 50 | 0.54% |
| Austria | 43 | 0.46% |
| Greece | 43 | 0.46% |
| Norway | 41 | 0.44% |
| United Arab Emirates | 40 | 0.43% |
| Ukraine | 38 | 0.41% |
| Morocco | 35 | 0.38% |
| Nigeria | 35 | 0.38% |
| Sweden | 33 | 0.35% |
| Azerbaijan | 31 | 0.33% |
| Tunisia | 29 | 0.31% |
| New Zealand | 26 | 0.28% |
| Serbia | 25 | 0.27% |
| Chile | 23 | 0.25% |
| Denmark | 23 | 0.25% |
| Ireland | 22 | 0.24% |
| Czech Republic | 21 | 0.23% |
| Indonesia | 20 | 0.21% |
| Oman | 20 | 0.21% |
| Yemen | 20 | 0.21% |
| Colombia | 19 | 0.20% |
| Bangladesh | 17 | 0.18% |
| Finland | 16 | 0.17% |
| Slovenia | 15 | 0.16% |
| Lebanon | 13 | 0.14% |
| Argentina | 11 | 0.12% |
| Estonia | 11 | 0.12% |
| Kazakhstan | 11 | 0.12% |
| Kuwait | 11 | 0.12% |
| Palestine | 10 | 0.11% |
| Qatar | 10 | 0.11% |
| Libya | 8 | 0.09% |
| Ethiopia | 9 | 0.10% |
| Cyprus | 7 | 0.08% |
| Macao | 6 | 0.06% |
| Slovakia | 6 | 0.06% |
| Croatia | 5 | 0.05% |
| Georgia | 5 | 0.05% |
| Cameroon | 4 | 0.04% |
| Ecuador | 4 | 0.04% |
| Lithuania | 4 | 0.04% |
| Tanzania | 4 | 0.04% |
| Trinidad and Tobago | 4 | 0.04% |
| Uruguay | 4 | 0.04% |
| Uzbekistan | 4 | 0.04% |
| Venezuela | 4 | 0.04% |
| Afghanistan | 3 | 0.03% |
| Bahrain | 3 | 0.03% |
| Dominican Republic | 3 | 0.03% |
| Peru | 3 | 0.03% |
| Philippines | 3 | 0.03% |
| Belarus | 2 | 0.02% |
| Cote d'Ivoire | 2 | 0.02% |
| Ghana | 2 | 0.02% |
| Kyrgyzstan | 2 | 0.02% |
| Malawi | 2 | 0.02% |
| North Korea | 2 | 0.02% |
| Yugoslavia | 2 | 0.02% |
| Albania | 1 | 0.01% |
| Armenia | 1 | 0.01% |
| Benin | 1 | 0.01% |
| Chad | 1 | 0.01% |
| Honduras | 1 | 0.01% |
| Jamaica | 1 | 0.01% |
| Kenya | 1 | 0.01% |
| Luxembourg | 1 | 0.01% |
| Madagascar | 1 | 0.01% |
| Mauritius | 1 | 0.01% |
| Moldova | 1 | 0.01% |
| Nepal | 1 | 0.01% |
| North Macedonia | 1 | 0.01% |
| Sierra Leone | 1 | 0.01% |
| Swaziland | 1 | 0.01% |
| Syrian Arab Republic | 1 | 0.01% |
| Undefined | 122 | 1.31% |

Two South Korean universities from Seoul top the list of authors' affiliations, followed by a Romanian, an Indian, a Chinese, a Saudi Arabian and a Turkish university (see Fig. 4).

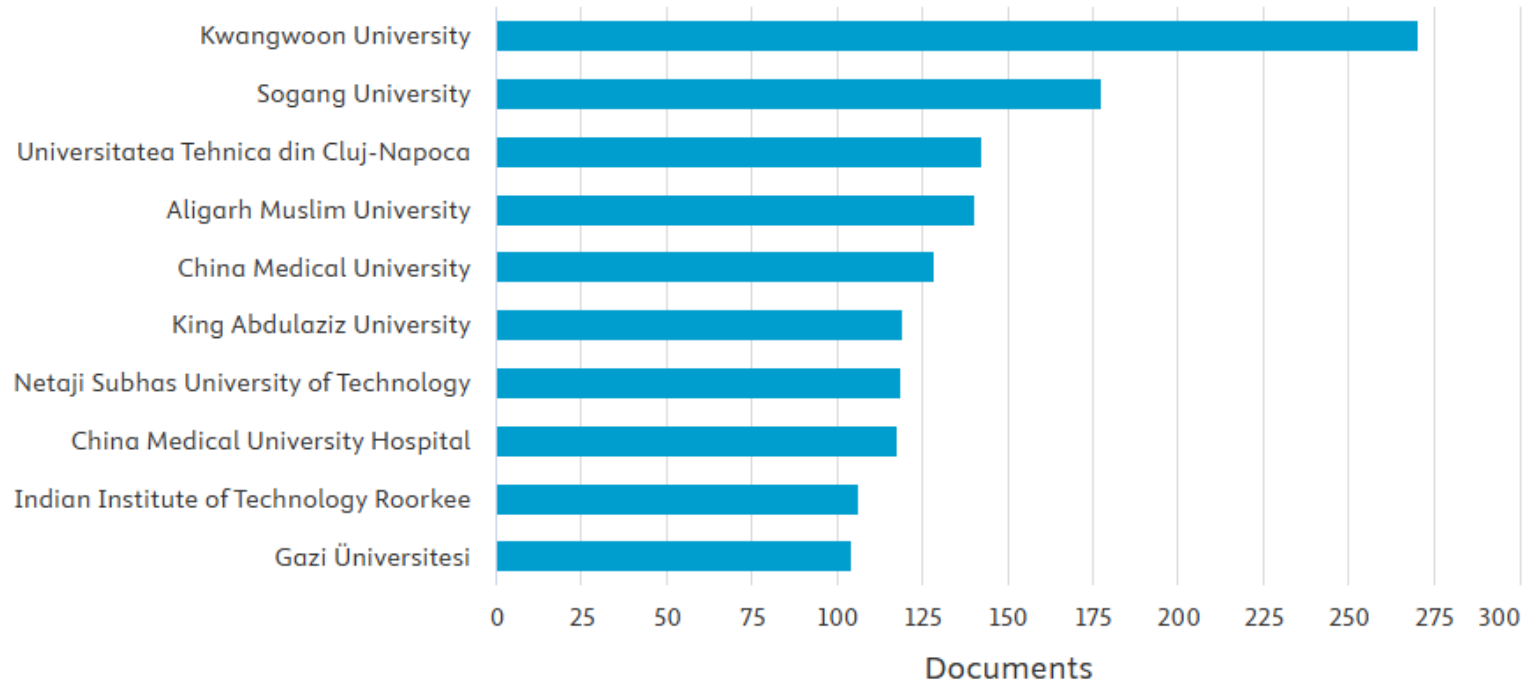

**Figure 4.** Top ten universities with documents containing the keyword "Bernstein polynomial" in all search fields.

Other universities and institutions mentioned in more than 30 documents include Ankara University (Turkey), Universitatea Babeș-Bolyai (Romania), Zhejiang University (China), Islamic Azad University (Iran), University of Oradea (Romania), Kyungpook National University (South Korea), Eastern Mediterranean University (Cyprus), Hannam University (South Korea), Akdeniz University (Turkey), French National Centre for Scientific Research (France), Shiraz University of Technology (Iran), University of Zaragoza (Spain), Çankaya University (Turkey), Ministry of Education of the People's Republic of China, Islamic Azad University, Karaj Branch (Iran), Shahid Beheshti University (Iran), Lucian Blaga University of Sibiu (Romania), Institute for Space Sciences, Bucharest (Romania), Pukyong National University (South Korea), Xiamen University (China), Malayer University (Iran), Atilim University (Turkey), Gaziantep University (Turkey), Bolu Abant İzzet Baysal University (Turkey).

Figure 5 shows the top ten academics and the co-authorship overlay visualisation for 2010-2025 (the colours change from blue for older documents to red for newer documents).

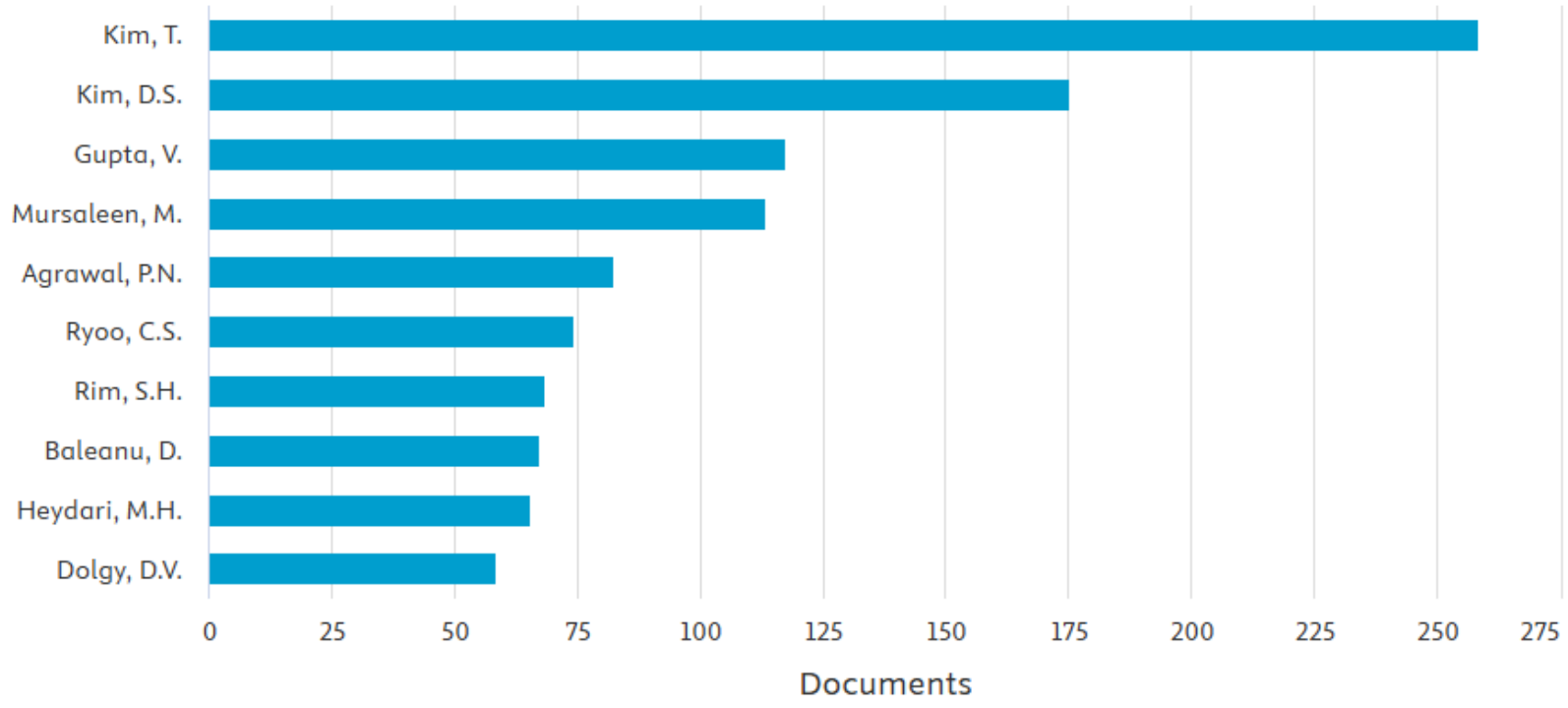

(**a**)

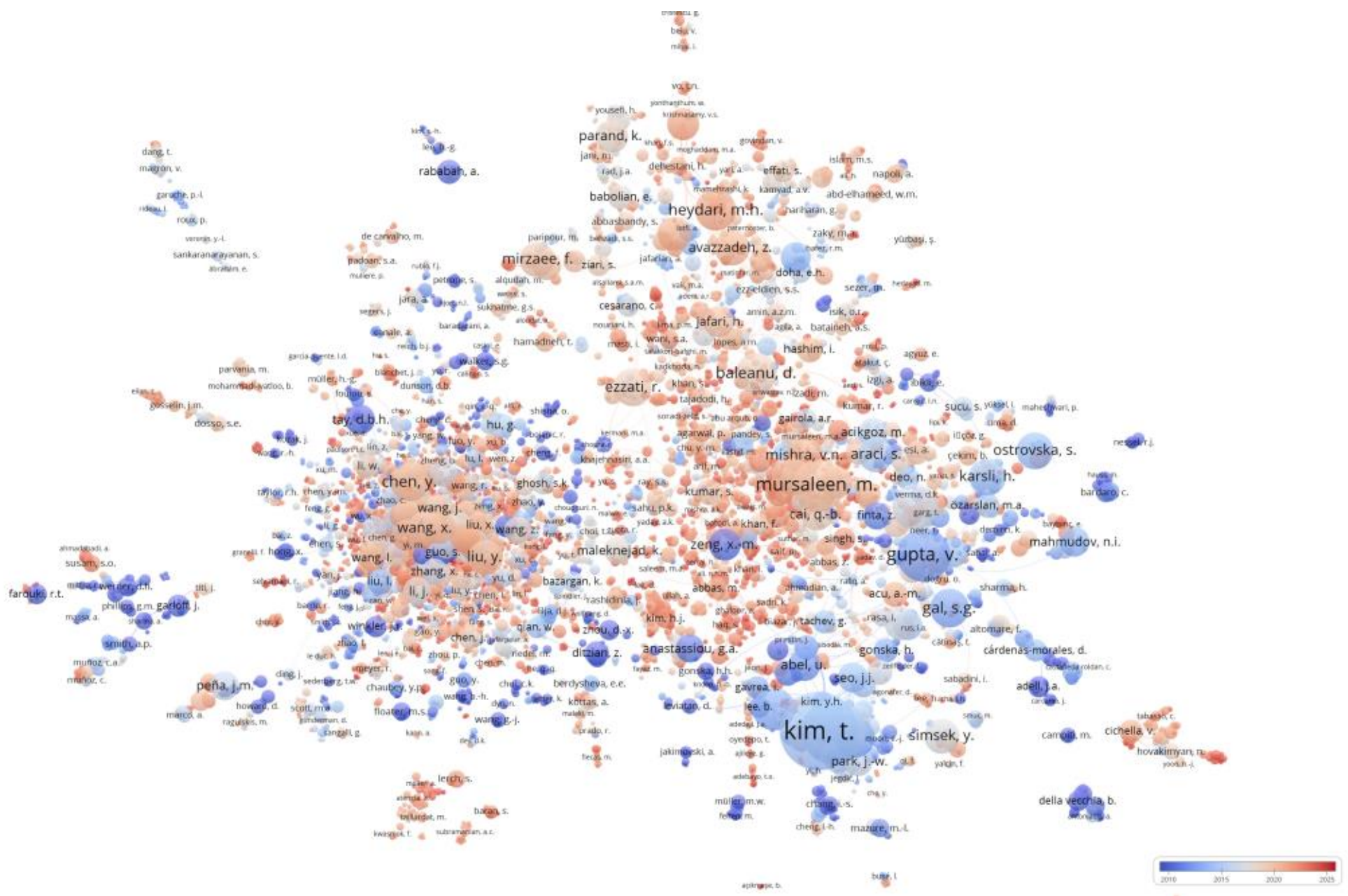

(**b**)

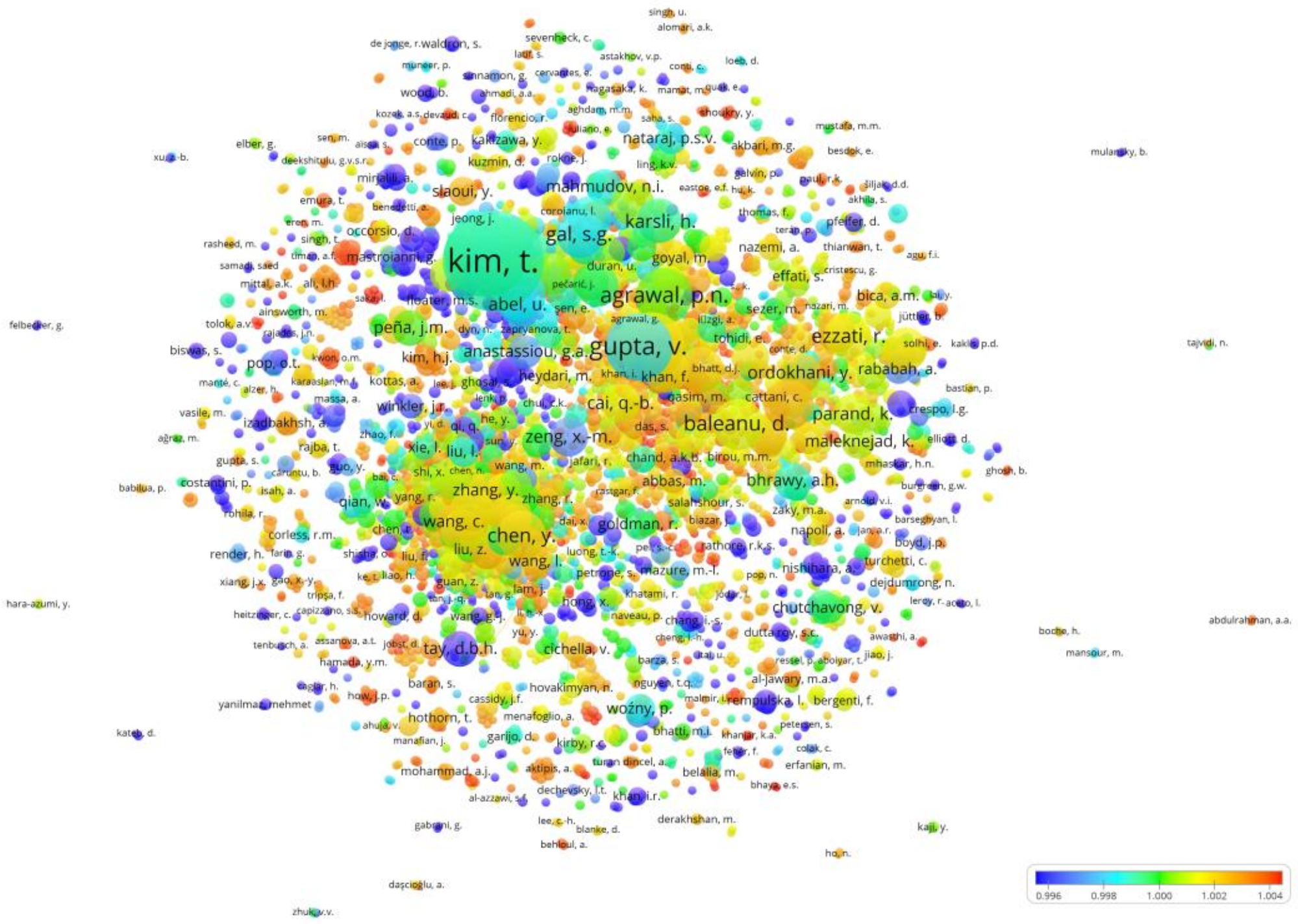

(**c**)

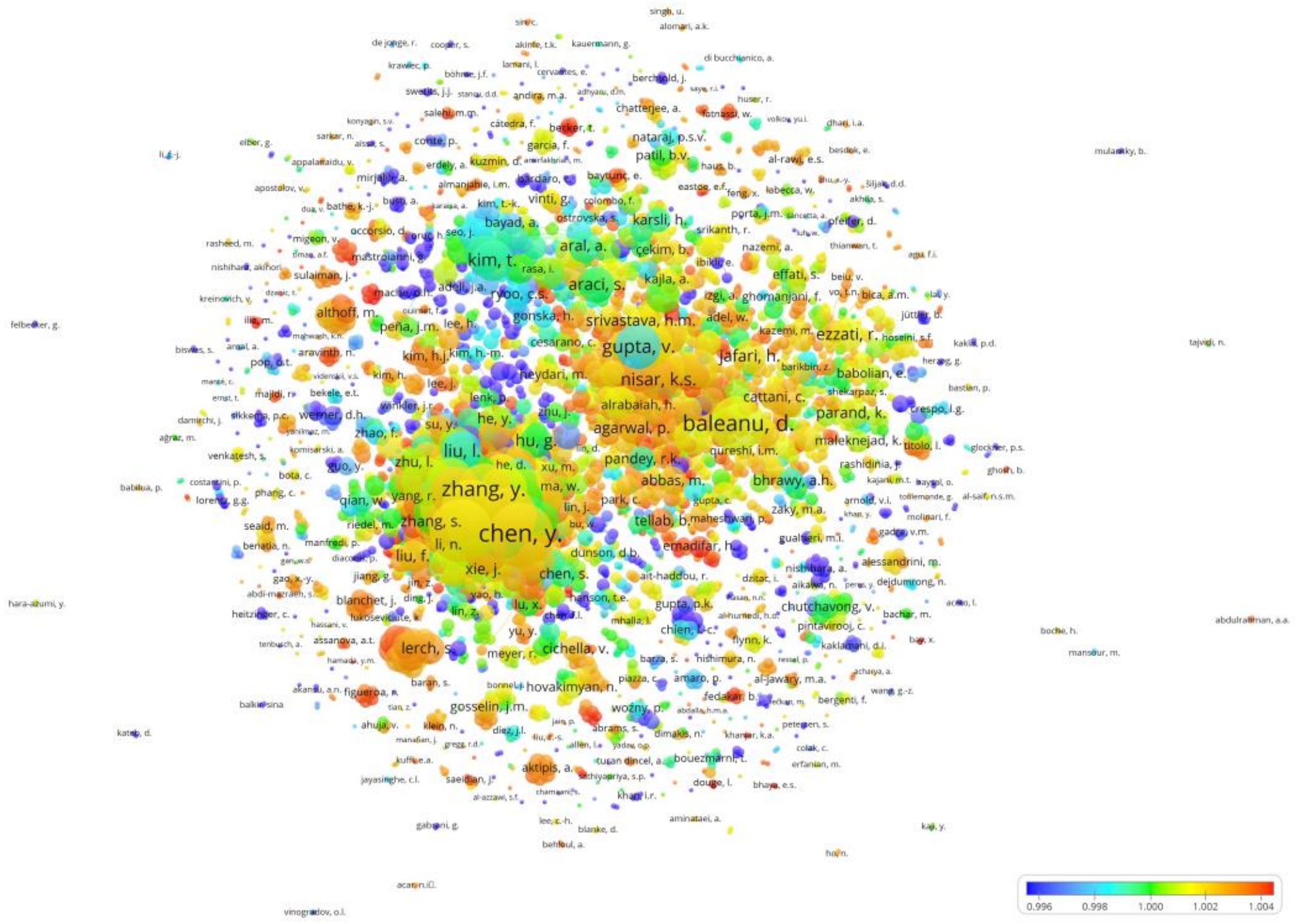

(**d**)

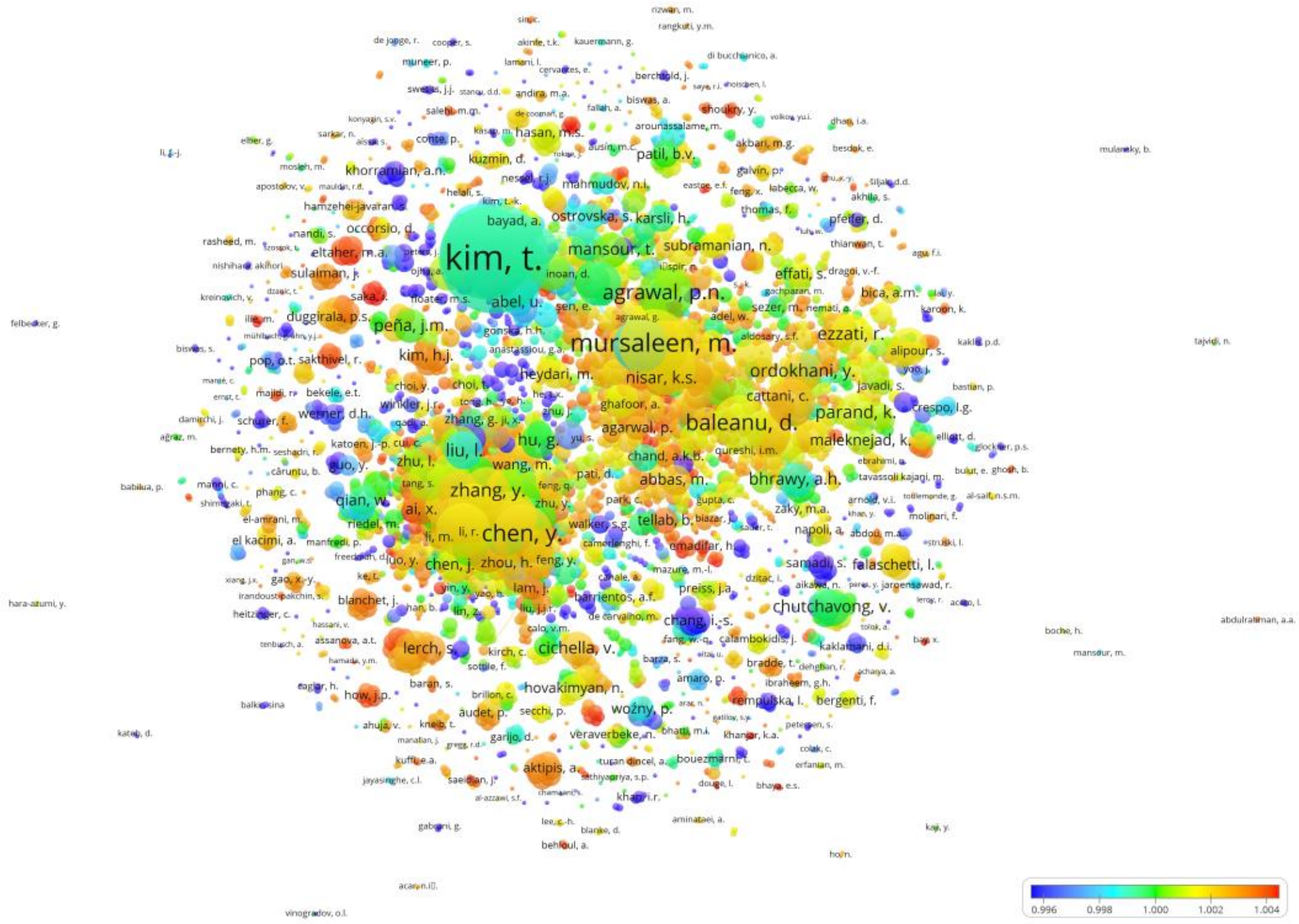

(**e**)

**Figure 5.** (a) Top ten academics whose documents contain the keyword "Bernstein polynomial" in all search fields, (b) Co-authorship overlay visualisation with a more dynamic variant available online at https://youtu.be/ipblMZ_saKA **(c)** Overlay visualisation in which the size of each node represents the weight based on the number of documents, while the colour, which changes from blue to red, reflects the normalised

average publication year[3], with blue representing older documents and red representing more recent ones. **(d)** Overlay visualisation in which the size of each node represents the weight based on the number of links, while the colour, which changes from blue to red, reflects the normalised average publication year, with blue representing older documents and red representing more recent ones. **(e)** Overlay visualisation in which the size of each node represents the weight based on the total link strength, while the colour, which changes from blue to red, reflects the normalised average publication year, with blue representing older documents and red representing more recent ones.

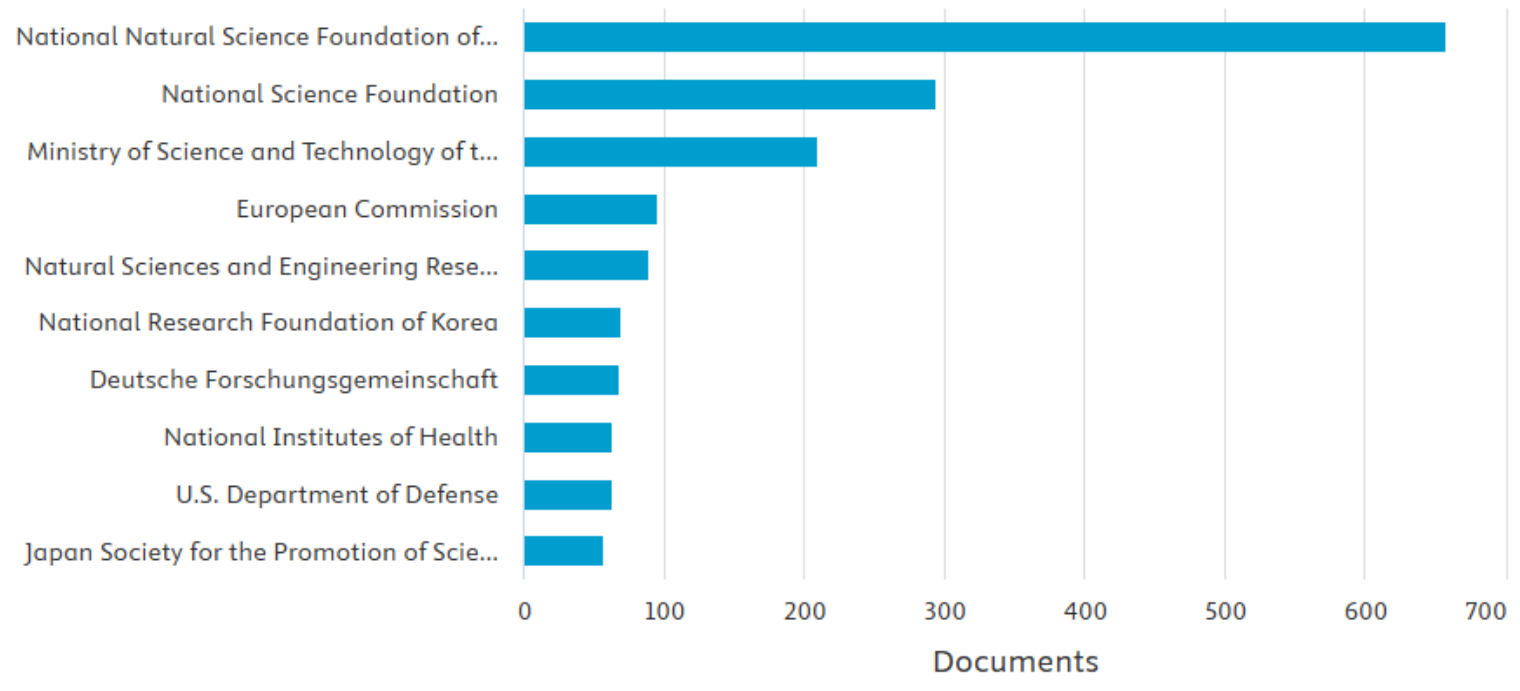

**Figure 6.** Top 10 funding agencies.

The list of funding agencies is shown in Figure 6 and Table 2. The leading funding agency is the National Natural Science Foundation of China, and the third one is the Ministry of Science and Technology of the People's Republic of China, which indicates the enormous support provided by China for research in areas related to Bernstein polynomials and their applications. Although the National Research Foundation of Korea is mentioned as a funding agency in only 68 documents, authors from this country are among the most active researchers mentioning the "Bernstein polynomial" in their documents (Fig. 5a).

**Table 2.** The list of funding sponsors mentioned in documents with the keywords "Bernstein polynomial" was found in any search field. Unfortunately, some funding agencies do not include the name of a country in their title, and some are written in a language other than English.

| Funding sponsor | Documents |
|---|---|
| National Natural Science Foundation of China | 657 |
| National Science Foundation | 293 |
| Ministry of Science and Technology of the People's Republic of China | 209 |
| European Commission | 94 |
| Natural Sciences and Engineering Research Council of Canada | 88 |
| National Research Foundation of Korea | 68 |
| Deutsche Forschungsgemeinschaft[4] | 67 |
| National Institutes of Health | 62 |
| U.S. Department of Defense | 62 |
| Japan Society for the Promotion of Science | 56 |
| Fundamental Research Funds for the Central Universities | 54 |
| University Grants Commission | 53 |
| Office of Naval Research | 49 |
| National Key Research and Development Program of China | 47 |
| Department of Science and Technology, Ministry of Science and Technology, India | 46 |
| Ministry of Education of the People's Republic of China | 42 |
| Natural Science Foundation of Fujian Province | 41 |

[3] According to the VOSviewer manual available online (https://www.vosviewer.com/documentation/Manual_VOSviewer_1.6.20.pdf), "Avg. pub. year" means the average publication year of documents containing a keyword or term, or the average publication year of documents published by a source, author, organisation or country.

[4] German Research Foundation

| | |
|---|---|
| Science and Engineering Research Board | 38 |
| European Regional Development Fund | 37 |
| Ministerio de Economía y Competitividad[5] | 37 |
| Air Force Office of Scientific Research | 36 |
| Council of Scientific and Industrial Research, India | 35 |
| U.S. Department of Energy | 35 |
| U.S. Department of Health and Human Services | 35 |
| Seventh Framework Programme | 34 |
| Engineering and Physical Sciences Research Council | 33 |
| Government of Canada | 33 |
| Natural Science Foundation of Zhejiang Province | 33 |
| Schweizerischer Nationalfonds zur Förderung der Wissenschaftlichen Forschung[6] | 33 |
| European Research Council | 32 |
| China Postdoctoral Science Foundation | 31 |
| Horizon 2020 Framework Programme | 31 |
| Ministry of Education, Culture, Sports, Science and Technology | 31 |
| Natural Science Foundation of Hebei Province | 31 |
| Ministry of Human Resource Development | 30 |
| Akdeniz University, Turkey | 29 |
| Directorate for Mathematical and Physical Sciences | 29 |
| U.S. Air Force | 28 |
| Air Force Materiel Command | 27 |
| Ministerio de Ciencia, Innovación y Universidades[7] | 26 |
| Ministero dell'Istruzione, dell'Università e della Ricerca[8] | 26 |
| U.S. Navy | 26 |
| Natural Science Foundation of Jiangsu Province | 24 |
| UK Research and Innovation | 23 |
| Istituto Nazionale di Alta Matematica "Francesco Severi"[9] | 22 |
| Junta de Andalucía[10] | 22 |
| National Aeronautics and Space Administration | 22 |
| Agencia Estatal de Investigación[11] | 21 |
| China Scholarship Council | 21 |
| Fundação para a Ciência e a Tecnologia[12] | 21 |

The second source of funding mentioned in 293 documents, the National Science Foundation (NSF), doesn't mention any country by name, but an analysis of the locations of the universities shows that 75.32% of all workplaces are in US universities, and 5.06% are Chinese, 2.74% are Canadian, 2.11% are German and British, 1.69% are French, and so on. Therefore, it's evident that the National Science Foundation (https://nsf.gov) is an independent agency of the United States federal government that supports fundamental research and education. However, for efficient tracking of research funding, the name of the fund mentioned in manuscripts should be unique, and various academic databases can also integrate these names into the document properties so that it would be easy to know that the development of a document was supported by research funding.

[5] Ministry of Economy and Competitiveness, Spain

[6] Swiss National Science Foundation, Switzerland

[7] Ministry of Science, Innovation and Universities, Spain

[8] Ministry of Education, University and Research, Italy

[9] The Francesco Severi National Institute of Higher Mathematics, Italy

[10] Regional Government of Andalusia, Spain

[11] State Research Agency, Spain

[12] The Foundation for Science and Technology, Portugal

Analysis of author names in NSF-supported documents using AI chatbots[13] revealed that some US-based academics have names of Chinese and other origins (Table 3), but different tools produced different results.

**Table 3.** AI-based analysis of author names using different LLMs performed on 23 January 2025.

| Name origin | Llama 3 | Perplexity | Mixtral | Average |
|---|---|---|---|---|
| United States | 34.86% | 35.00% | 30.60% | 33.49% |
| China | 20.41% | 20.00% | 14.60% | 18.34% |
| United Kingdom | 5.36% | 8.00% | 2.40% | 5.25% |
| Germany | 3.59% | 5.00% | 6.10% | 4.90% |
| Canada | 3.59% | 5.00% | 7.90% | 5.50% |
| France | 3.59% | 4.00% | 4.30% | 3.96% |
| India | 3.59% | 3.00% | 5.50% | 4.03% |
| Japan | 1.79% | 2.00% | 0.60% | 1.46% |
| Others | 23.22% | 18.00% | 28.00% | 23.07% |
| Total | 100.00% | 100.00% | 100.00% | |

Table 4 shows the top 10 most cited papers, including citation counts, journal names, and publishers. 40% of these top papers were published in Computer Aided Geometric Design, a key Elsevier journal for research in the mathematical foundations of free-form curves, surfaces and solids. A visualisation of the evolution of sources that published ten or more documents related to the keyword "Bernstein polynomial" in data extracted from Scopus can be found online at https://youtu.be/CTZv1iTT7gA.

**Table 4.** The top 10 most cited articles that contain the keyword "Bernstein polynomial" in the article title, abstract or keywords according to the Scopus database. The article by Sederberg & Parry is indexed twice in Scopus, showing a different number of citations.

| Authors | Year | Title | Source title | Cited by | Publisher |
|---|---|---|---|---|---|
| Sederberg Thomas W.; Parry Scott R. [7] | 1986 | Free-form deformation of solid geometric models | Computer Graphics (ACM) | 2285 | Association for Computing Machinery |
| Farin G. [8] | 1986 | Triangular Bernstein-Bézier patches | Computer Aided Geometric Design | 497 | Elsevier |
| Farouki R.T. [9] | 2012 | The Bernstein polynomial basis: A centennial retrospective | Computer Aided Geometric Design | 379 | Elsevier |
| Farouki R.T.; Rajan V.T. [10] | 1988 | Algorithms for polynomials in Bernstein form | Computer Aided Geometric Design | 257 | Elsevier |
| Farouki R.T.; Rajan V.T. [11] | 1987 | On the numerical condition of polynomials in Bernstein form | Computer Aided Geometric Design | 241 | Elsevier |
| King J.P. [12] | 2003 | Positive linear operators which preserve $x^2$ | Acta Mathematica Hungarica | 236 | Hungarian Academy of Sciences |
| Gao F.; Wu W.; Lin Y.; Shen S. [13] | 2018 | Online Safe Trajectory Generation for Quadrotors Using Fast Marching Method and Bernstein Basis Polynomial | Proceedings - IEEE International Conference on Robotics and Automation | 212 | Institute of Electrical and Electronics Engineers Inc. |
| Sederberg T.W.; Parry S.R. [14] | 1986 | Free-form deformation of solid geometric models | Proceedings of the 13th Annual Conference on Computer Graphics and Interactive Techniques, SIGGRAPH 1986 | 193 | Association for Computing Machinery |
| Bhrawy A.H.; Taha T.M.; Machado J.A.T. [15] | 2015 | A review of operational matrices and spectral techniques for fractional calculus | Nonlinear Dynamics | 181 | Kluwer Academic Publishers |

[13] It is not possible to use AI chatbots for all 9310 documents, as most chatbots are not able to process very long texts.

| | | | | | |
|---|---|---|---|---|---|
| Ostrovska S. [16] | 2003 | $q$-Bernstein polynomials and their iterates | Journal of Approximation Theory | 180 | Academic Press Inc. |

## 3. Network Visualisations and Analysis

Network visualisations and analysis were performed using VOSviewer (https://www.vosviewer.com/), a freeware tool for constructing and visualising bibliometric networks. Gephi (https://gephi.org/) is another free and open-source tool, the leading visualisation and exploration software for all graphs and networks. The minimum number of occurrences of a keyword was set to 3, and 5087 keywords met the threshold. VOSviewer has created 21 clusters that can be seen in different colours in Fig. 7a, 161,307 links, with a total link strength of 261,230. Three types of network visualisations are shown in Fig. 7.

In VOSviewer [17-19], we used the following metrics:

- ***Links*** refer to the number of direct connections or relationships between nodes in the graph. In social networks, for example, this is like counting the number of friendships between individuals, but without considering the strength of those friendships.

- ***Total link strength*** is the cumulative weight of all links between nodes in the graph. The weights can represent the intensity, frequency or importance of the relationships. For example, a citation network could reflect the number of times two articles have been co-cited.

- ***Occurrences*** refer to the number of times a particular keyword appears in the dataset used to construct the graph. For example, the occurrence of "Bernstein polynomial" indicates how often this term is mentioned in the data, regardless of its connections to other terms.

Figure 7 contains network visualisations for the keywords studied. Fig. 7a is a basic network visualisation where we can observe several clusters coloured by different colours, Fig. 7b is a density visualisation (or density map) where the red colour shows the most popular topics based on the number of occurrences, and Fig. 7c is the most interesting because the years from 1949 to 2025 were normalised and connected to the colour map used, so we can observe that keywords in blue colours, such as algorithms, approximation theory, filters, CAD, convexity, image reconstruction, and others were relatively more common for older manuscripts, and motion planning, trajectories, uncertainty, machine learning, weather forecasting, collision avoidance, neural networks, continuous-time systems, and others can be found in very recent documents. This means a shift from traditional mathematical topics to applications and numerical methods.

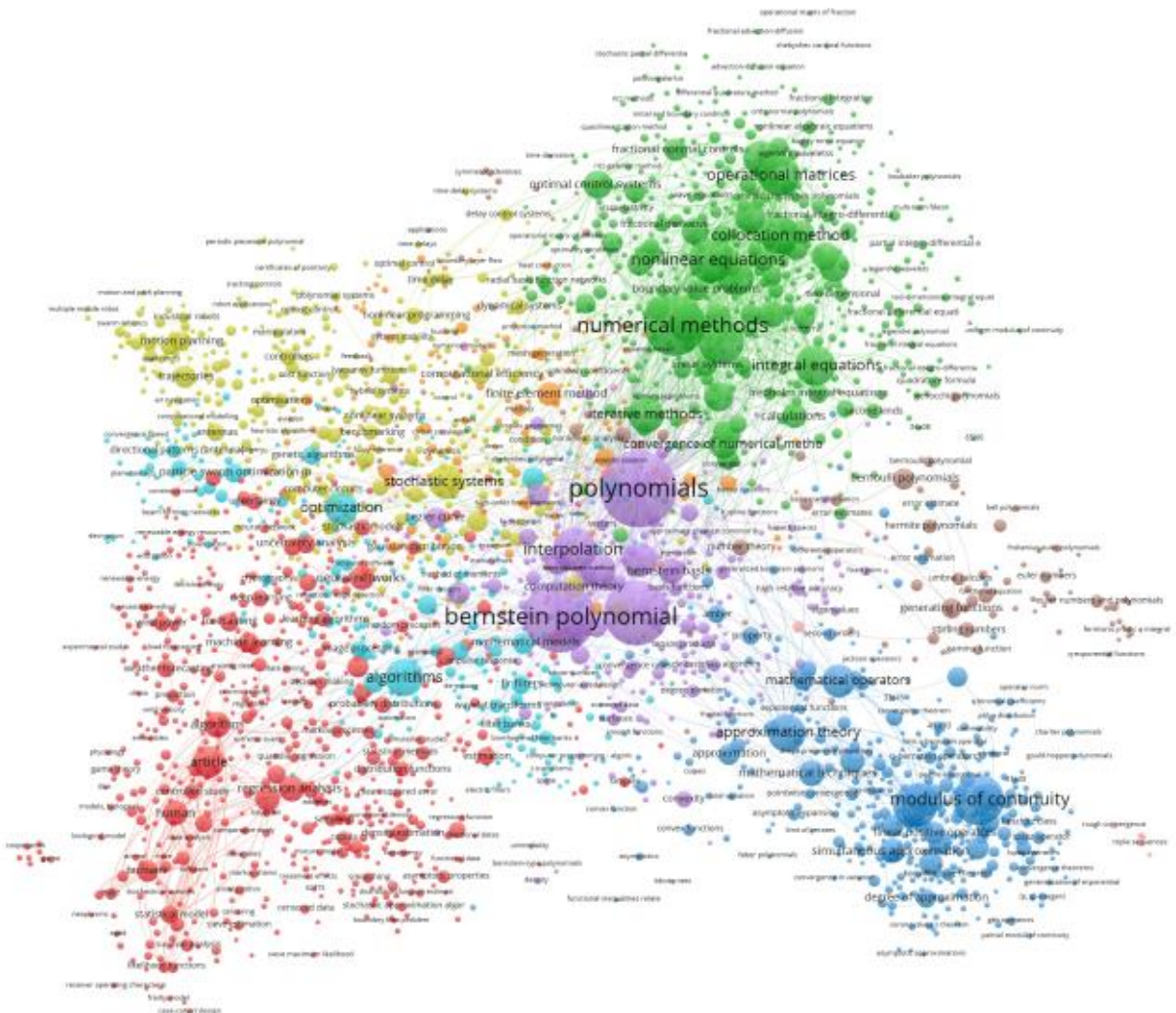

(**a**)

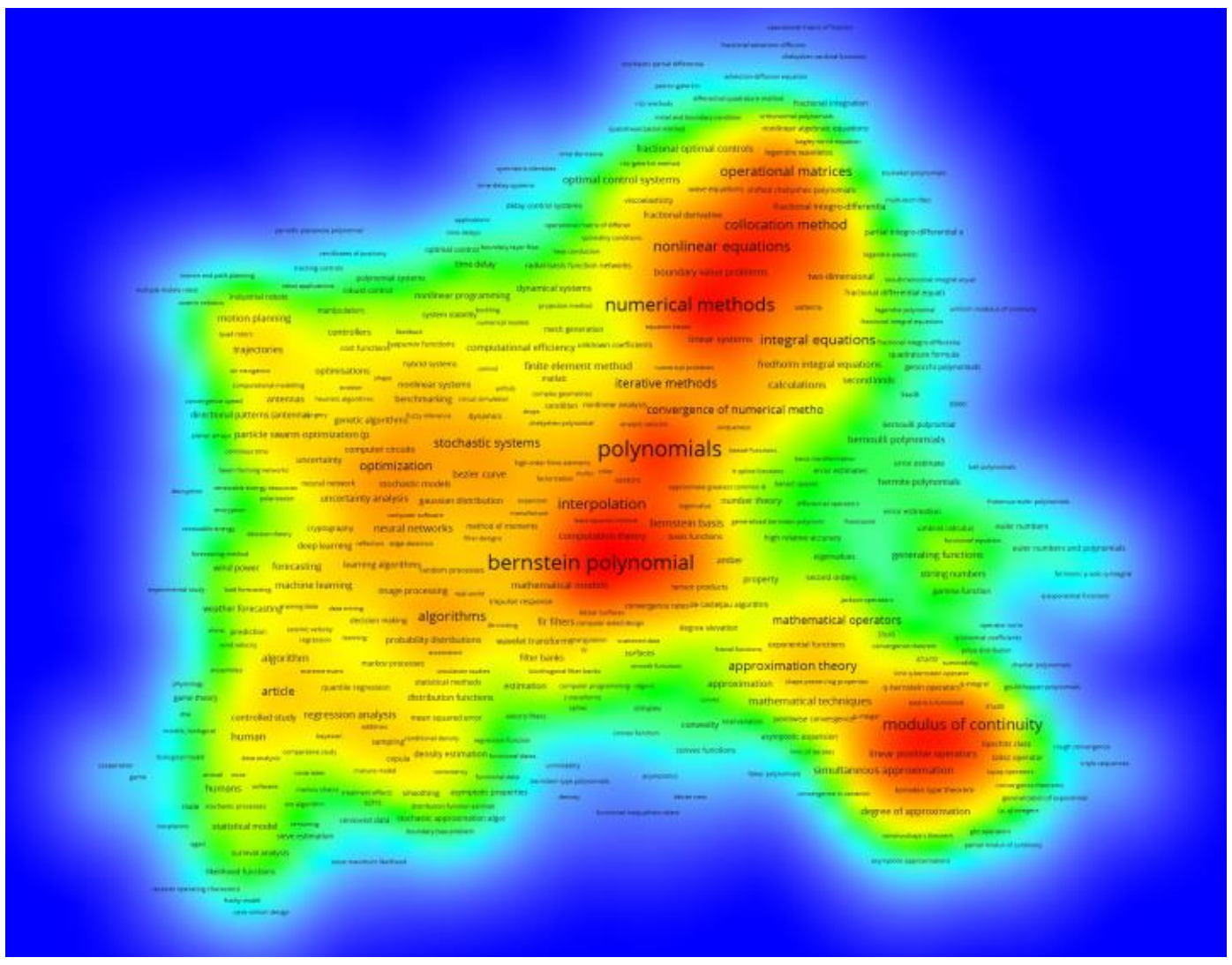

(**b**)

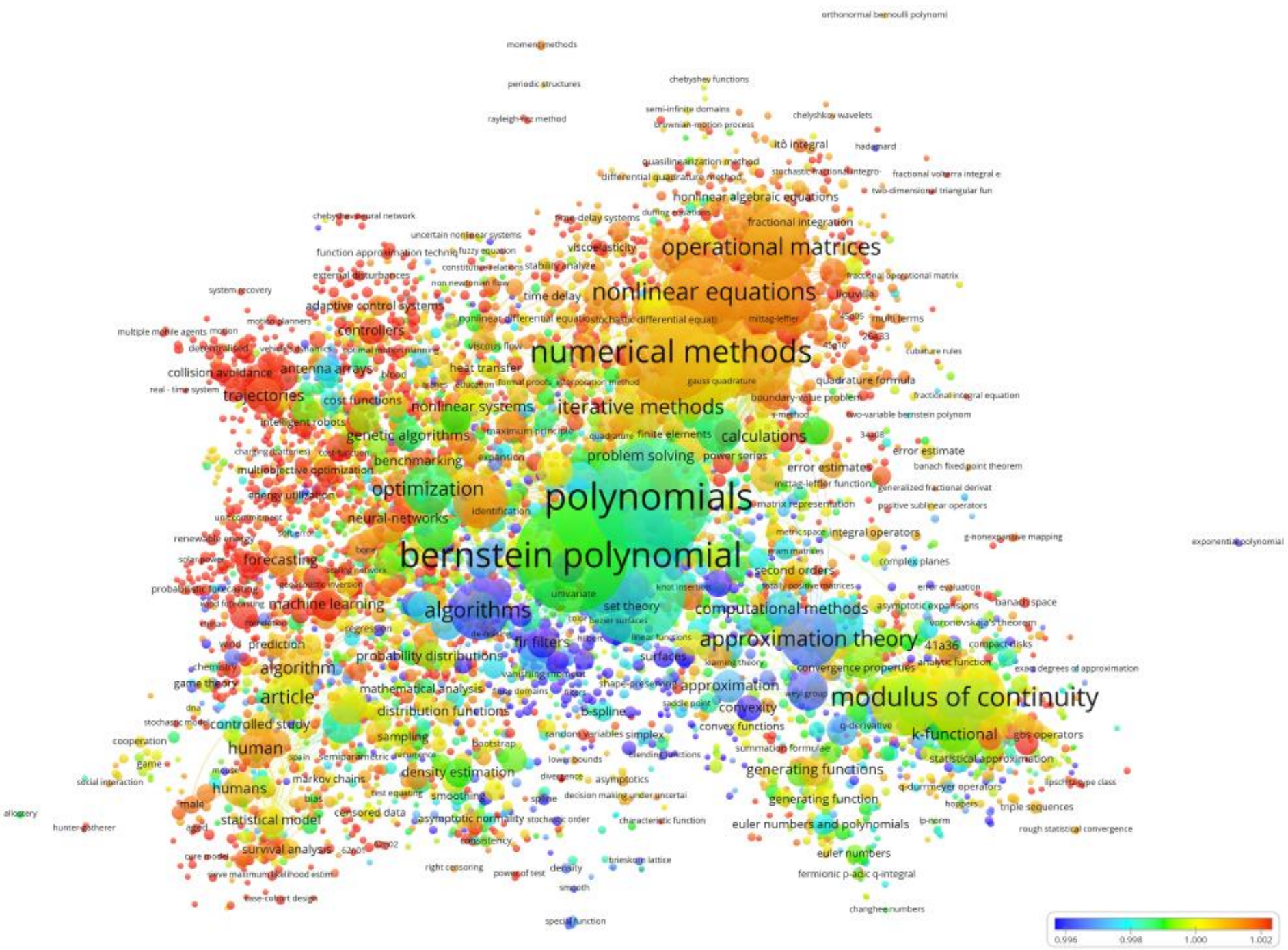

(**c**)

**Figure 7.** Network visualisation based on the Scopus query and 9310 manuscripts found as of 23 January 2025. **(a)** Network visualisation, **(b)** Density visualisation. A more dynamic visualisation with zooming is available at https://youtu.be/bpqYJVEynQk **(c)** Overlay visualisations where years from 1949 to 2025 are normalised by dividing by the mean. Here, cool colours represent older research topics, while warm colours represent more recent research topics, so the visualisation helps to understand the trends for the related research keywords.

In network visualisation (Fig. 8), there are two phrases related to Bernstein polynomials: singular and plural. Both show stronger links with topics such as trajectories, prediction, optimisation, algorithms, approximation theory, modulus of continuity, mathematical operators, numerical and iterative methods, matrix algebra, operational matrices, etc. For "Bernstein polynomial" (Fig. 8b): Links=2561, Total link strength=7911, Occurrences=888, and for "Bernstein polynomials" (Fig. 8a): Links=1784, Total link strength=4717, Occurrences=728.

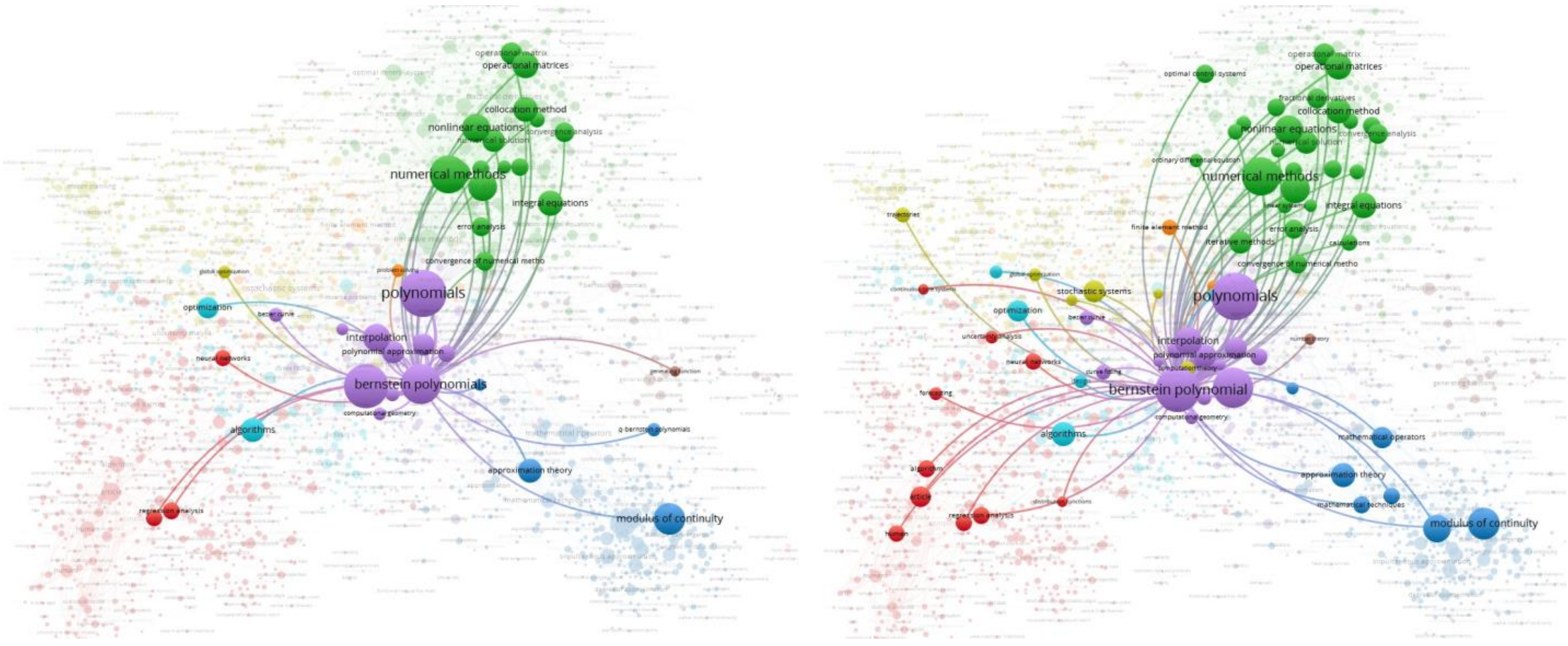

(**a**) (**b**)

**Figure 8.** Links of (a) "Bernstein polynomials" and (b) "Bernstein polynomials" in the resulting network visualisation.

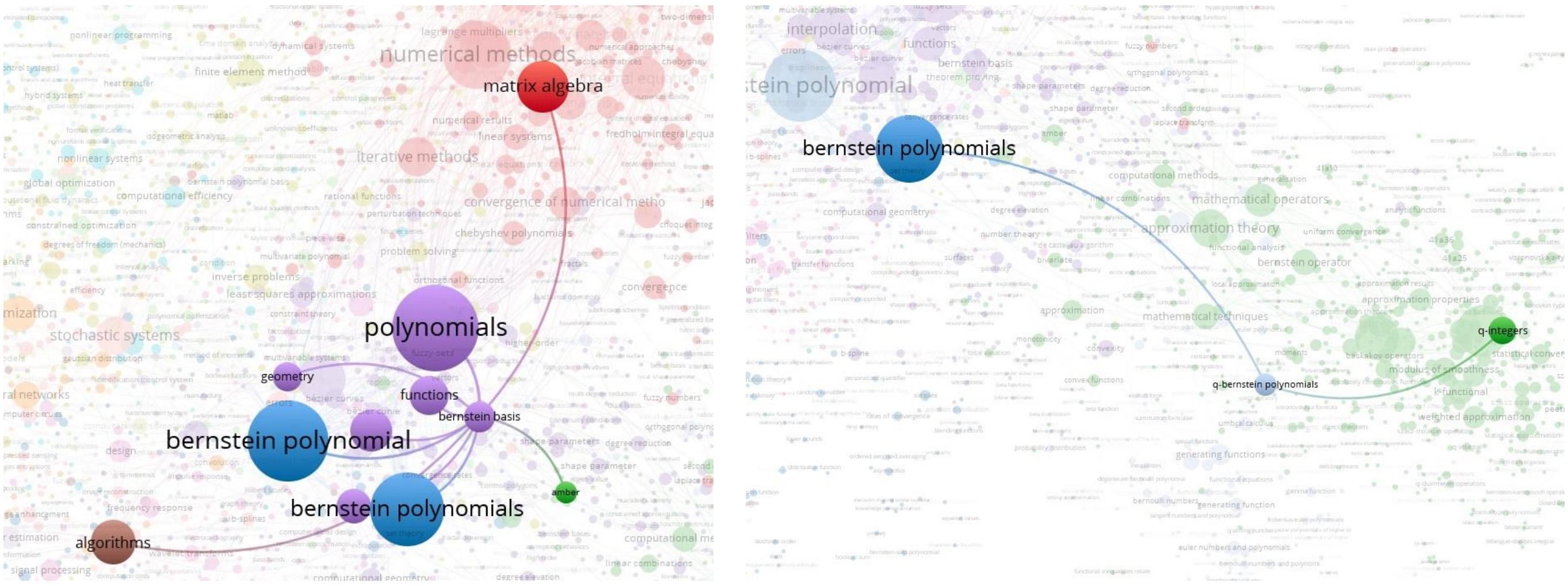

(**a**) (**b**)

**Figure 9.** Links of (a) "Bernstein basis" (Links=504, Total link strength=1081) and (b) "$q$-Bernstein polynomial" (Links=162, Total link strength=315) in the resulting network visualisation.

A-K-J Sankey diagram, which visualises the main authors (A), keywords (K), journals (J) and how they are related, can be found at https://youtu.be/d2hnhklw8io

## 4. Areas Where Bernstein Polynomials Are Rarely Used

By searching for different keywords in VOSviewer, we can find areas not well connected with Bernstein polynomials in documents from the Scopus database. For example, searching for "biology" returned two items from two clusters (Fig. 10).

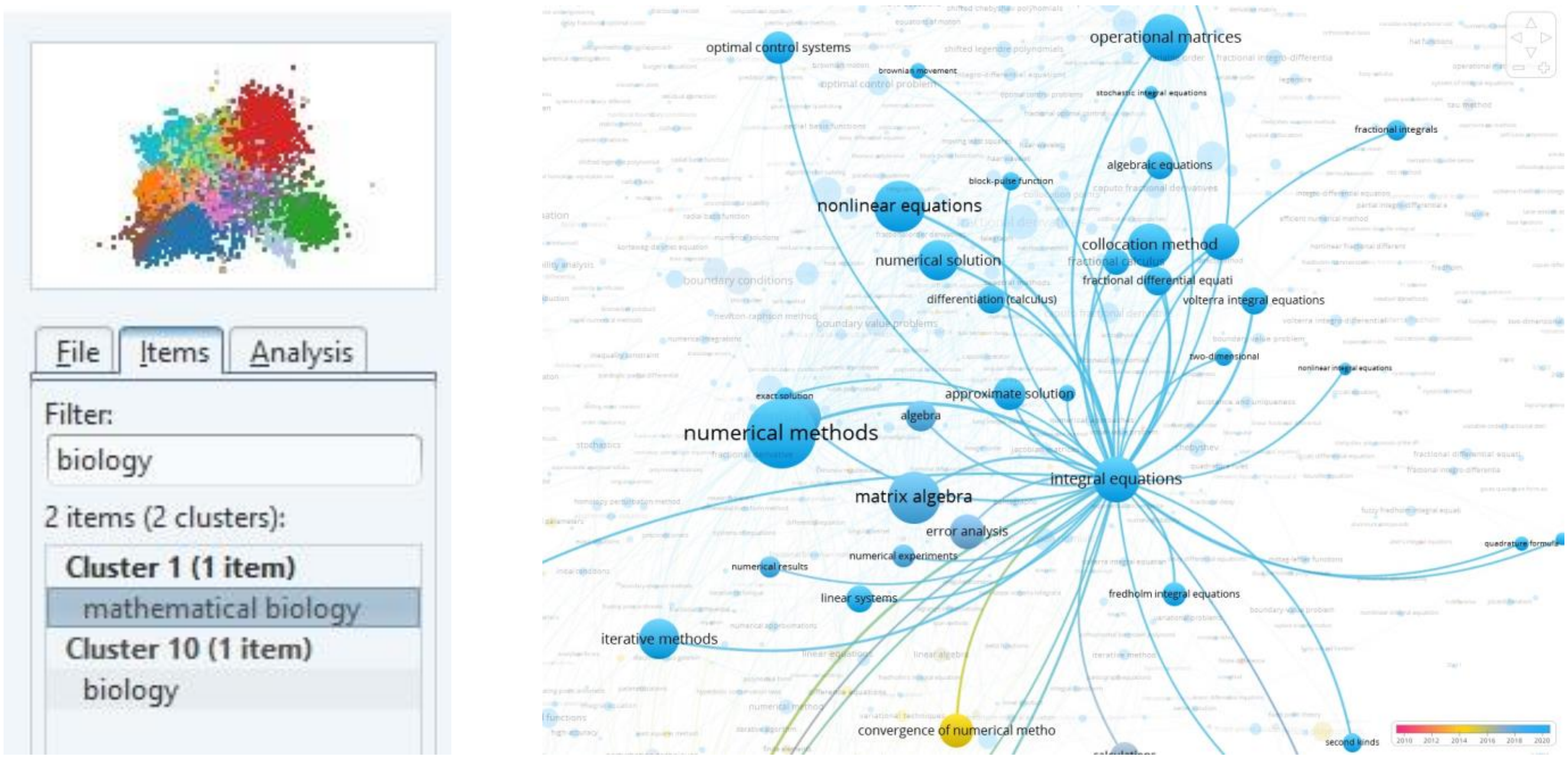

(**a**) (**b**)

**Figure 10.** (a) One of the interesting branches of applied mathematics, "mathematical biology" (Links=35, Total link strength=41, Occurrences=3), doesn't have a strong connection with "Bernstein polynomial". We have something similar for "biology" Links=38, Total link strength=39, Occurrences=3. (b) The item "mathematical biology" has a relatively smaller scale than "integral equations", a mathematical tool used to construct mathematical models, particularly in biology and ecology.

A Scopus search for (TITLE-ABS-KEY("Bernstein polynomial") AND TITLE-ABS-KEY("mathematical biology")) returned only four articles, of which only two papers [20, 21] have titles related to mathematical biology. Examples of other keywords are given in Table 5. Strangely, not only "Bernstein polynomial" but both "Bézier curve" and "Bézier surface" appear together with "aerospace engineering" in only fifteen articles. However, the numbers can increase significantly when searching in all fields. A Scopus search in all fields (ALL("Bernstein polynomial") AND ALL("prize" OR "award")) returned only ten documents, most of which have nothing to do with mathematical prizes or awards.

**Table 5.** Keywords rarely appear with "Bernstein polynomial" in document titles, abstracts or keywords in the Scopus database.

| Keyword searched with "Bernstein polynomial" in Scopus in document title, abstract or keywords[14] | Number of documents found | Notes |
|---|---|---|
| robotics | 21 | Only nine articles contain the word "robot" in their title |
| aerospace engineering | 4 | Only two articles are related to aerodynamics and flight vehicle |
| composite material | 4 | Only two articles mention composites in the title |
| fuzzy logic | 6 | Four articles contain the word "fuzzy" in the title |
| electrical engineering | 3 | Only one work belongs to a mathematical journal |
| energy engineering | 0 | No documents were found |
| nuclear | 3 | Two documents published in bioscience sources |
| mechanical engineering | 1 | Recent work published in 2023 |
| special function[15] | 9 | Four articles were published in 2024 |
| Fields medal | 0 | No documents were found |
| Shaw prize | 0 | No documents were found |
| Abel prize | 0 | No documents were found |
| "Lobachevsky Prize" or "Lobachevsky Medal" | 0 | No documents were found |
| Turing Award | 0 | No documents were found |
| Nobel prize | 0 | No documents were found |

I've started a chat with AI chatbots about possible applications of Bernstein polynomials in the fields mentioned in Tables 5 and 7. Interested readers can check out the links in Table 6.

**Table 6.** Suggestions from various AI chatbots about applications of Bernstein polynomials in different fields.

| AI chatbot | Responses |
|---|---|
| ChatGPT | https://bit.ly/3EjGUlu |
| Google Gemini | https://bit.ly/3PXoHtg |

14 Academics themselves can try to find such a field or a topic not well related to Bernstein polynomials by using a Scopus query (TITLE-ABS-KEY("Bernstein polynomial") AND TITLE-ABS-KEY("topic of your interest")), or by using "Bernstein polynomial" AND "topic of your interest" in Scholar Google.

15 This smaller number of documents doesn't mean that Bernstein polynomials or their generalisations cannot be represented as a special case of more general special functions or in terms of other polynomials, such as Jacobi polynomials [22].

| | |
|---|---|
| Perplexity | https://bit.ly/3Ei5NoQ |
| Microsoft Copilot | https://bit.ly/42B1R5u |

More detailed Table 7 shows the co-occurrence of a given keyword related to the Mathematical Subjects Classification and the keyword "Bernstein polynomial" in documents in the Scopus database. These numbers have been manually extracted from the Scopus database. The information can be helpful for academics looking for new or underexplored areas where Bernstein polynomials can be applied. Interestingly, the Pearson correlation coefficient between the number of documents found by searching in "All fields" and "Article title, abstract, keywords" equals 0.848, indicating a powerful and positive relationship between the two variables.

**Table 7.** Co-occurrence of a given keyword(s) related to the Mathematical Subjects Classification and "Bernstein polynomial" in documents in the Scopus database. The two columns "All fields" and "Article title, abstract, keywords" use the same colour map but are applied separately to each other, where red colour indicates low values, and green indicates high values (number of documents).

| Top-level code | Mathematics Subject Classification (MSC) and used keywords (in " ") | All fields | Article title, abstract, keywords |
|---|---|---|---|
| 00 | General | | |
| | "recreational mathematics" | 0 | 0 |
| | "philosophy of mathematics" | 0 | 0 |
| | "mathematical modelling" OR "mathematical modelling" | 1291 | 9 |
| 01 | History and biography | | |
| | "history" and "biography" | 7 | 0 |
| 03 | Mathematical logic and foundations | | |
| | "mathematical logic" | 4 | 0 |
| 05 | Combinatorics | | |
| | "combinatorics" | 248 | 2 |
| 06 | Order, lattices, ordered algebraic structures | | |
| | "lattice" | 277 | 12 |
| | "ordered algebraic structure" | 0 | 0 |
| 08 | General algebraic systems | | |
| | "algebraic system" | 115 | 16 |
| 11 | Number theory | | |
| | "number theory" | 383 | 25 |
| 12 | Field theory and polynomials | | |
| | "field theory" | 60 | 0 |
| 13 | Commutative algebra (Commutative rings and algebras) | | |
| | "commutative algebra" | 20 | 0 |
| | "commutative ring" | 11 | 0 |
| 14 | Algebraic geometry | | |
| | "algebraic geometry" | 108 | 1 |
| 15 | Linear and multilinear algebra; matrix theory | | |
| | "linear algebra" | 348 | 13 |
| | "multilinear algebra" | 35 | 0 |
| | "matrix theory" | 29 | 0 |
| 16 | Associative rings and (associative) algebras | | |
| | "associative ring" | 0 | 0 |
| | "associative algebra" | 1 | 0 |
| 17 | Non-associative rings and (non-associative) algebras | | |
| | "non-associative ring" | 0 | 0 |
| | "non-associative algebra" | 0 | 0 |
| 18 | Category theory; homological algebra | | |
| | "category theory" | 2 | 0 |
| | "homological algebra" | 4 | 0 |
| 19 | K-theory | | |
| | "K-theory" | 4 | 0 |
| 20 | Group theory and generalisations | | |
| | "group theory" | 33 | 0 |
| 22 | Topological groups, Lie groups (and analysis upon them) | | |
| | "topological group" | 2 | 0 |
| | "Lie group" | 65 | 1 |
| 26 | Real functions (including derivatives and integrals) | | |
| | "real function" | 77 | 8 |
| 28 | Measure and integration | | |
| | "integration" | 1703 | 107 |
| 30 | Functions of a complex variable (including approximation theory in the complex domain) | | |
| | "complex variable" | 114 | 2 |
| | "complex domain" | 79 | 2 |
| | "approximation theory" | 2788 | 98 |
| 31 | Potential theory | | |
| | "potential theory" | 51 | 1 |
| 32 | Several complex variables and analytic spaces | | |
| | "several complex variables" | 15 | 0 |
| | "analytic space" | 9 | 0 |
| 33 | Special functions | | |
| | "special function" | 776 | 9 |
| 34 | Ordinary differential equations | | |
| | "ordinary differential equation" | 539 | 39 |
| | "ODE" | 219 | 15 |
| 35 | Partial differential equations | | |
| | "partial differential equation" | 1351 | 62 |
| | "PDE" | 137 | 8 |
| 37 | Dynamical systems and ergodic theory | | |
| | "dynamical system" | 558 | 17 |
| | "ergodic theory" | 13 | 1 |
| 39 | Difference (equations) and functional equations | | |
| | "difference equation" | 1227 | 5 |
| | "functional equation" | 182 | 20 |
| 40 | Sequences, series, summability | | |
| | "sequence" | 1668 | 161 |
| | "series" | 4080 | 130 |
| | "summability" | 280 | 18 |
| 41 | Approximations and expansions | | |
| | "approximation" | 6026 | 765 |
| | "expansion" | 1268 | 107 |
| 42 | Harmonic analysis on Euclidean spaces (including Fourier analysis, Fourier transforms, trigonometric approximation, trigonometric interpolation, and orthogonal functions) | | |
| | "harmonic analysis" | 186 | 3 |
| | "Fourier analysis" | 241 | 5 |
| | "Fourier transform" | 187 | 9 |
| | "trigonometric approximation" | 54 | 2 |
| | "trigonometric interpolation" | 16 | 2 |
| | "orthogonal function" | 194 | 14 |
| 43 | Abstract harmonic analysis | | |
| | "abstract harmonic analysis" | 1 | 0 |
| 44 | Integral transforms, operational calculus | | |
| | "integral transform" | 370 | 1 |
| | "operational calculus" | 30 | 0 |
| 45 | Integral equations | | |
| | "integral equation" | 1473 | 138 |
| 46 | Functional analysis (including infinite-dimensional holomorphy, integral transforms in distribution spaces) | | |
| | "functional analysis" | 1017 | 14 |
| | "infinite-dimensional holomorphy" | 0 | 0 |
| | "integral transform in distribution space" | 0 | 0 |
| 47 | Operator theory | | |
| | "operator theory" | 552 | 5 |
| 49 | Calculus of variations and optimal control; optimisation (including geometric integration theory) | | |
| | "calculus of variations" | 201 | 4 |
| | "optimal control" | 812 | 48 |
| | "optimization" or "optimisation" | 2585 | 237 |
| | "geometric integration theory" | 3 | 0 |
| 51 | Geometry | | |
| | "geometry" | 958 | 124 |
| 52 | Convex (geometry) and discrete geometry | | |
| | "convex geometry" | 6 | 0 |
| | "discrete geometry" | 3 | 0 |
| 53 | Differential geometry | | |
| | "differential geometry" | 79 | 2 |
| 54 | General topology | | |
| | "general topology" | 22 | 0 |
| 55 | Algebraic topology | | |
| | "algebraic topology" | 7 | 0 |
| 57 | Manifolds and cell complexes | | |
| | "manifold" | 201 | 9 |
| | "cell complex" | 1 | 0 |
| 58 | Global analysis, analysis on manifolds (including infinite-dimensional holomorphy) | | |
| | "global analysis" | 27 | 0 |
| | "analysis on manifolds" | 8 | 0 |
| 60 | Probability theory and stochastic processes | | |
| | "probability theory" | 313 | 8 |
| | "stochastic process" | 284 | 11 |
| 62 | Statistics | | |
| | "statistics" | 2223 | 89 |
| 65 | Numerical analysis | | |
| | "numerical analysis" | 1465 | 23 |
| 68 | Computer science | | |
| | "computer science" | 1841 | 6 |
| 70 | Mechanics of particles and systems (including particle mechanics) | | |
| | "mechanics of particles" | 0 | 0 |
| | "mechanics of systems" | 0 | 0 |
| 74 | Mechanics of deformable solids | | |
| | "mechanics of deformable solids" | 0 | 0 |
| 76 | Fluid mechanics | | |
| | "fluid mechanics" | 270 | 7 |
| 78 | Optics, electromagnetic theory | | |
| | "optics" | 204 | 3 |
| | "electromagnetic theory" | 29 | 0 |
| 80 | Classical thermodynamics, heat transfer | | |
| | "thermodynamics" | 50 | 1 |
| | "classical thermodynamics" | 0 | 0 |
| | "heat transfer" | 348 | 14 |
| 81 | Quantum theory | | |
| | "quantum theory" | 25 | 3 |
| 82 | Statistical mechanics, structure of matter | | |
| | "statistical mechanics" | 432 | 0 |
| | "structure of matter" | 0 | 0 |
| 83 | Relativity and gravitational theory (including relativistic mechanics) | | |
| | "relativity theory" or "theory of relativity" | 2 | 0 |
| | "gravitational theory", "theory of gravity", or "theory of gravitation" | 2 | 0 |
| | "relativistic mechanics" | 0 | 0 |
| 85 | Astronomy and astrophysics | | |
| | "astronomy" | 163 | 1 |
| | "astrophysics" | 139 | 3 |
| 86 | Geophysics | | |
| | "geophysics" | 107 | 1 |
| 90 | Operations research, mathematical programming | | |
| | "operations research" | 199 | 0 |
| | "mathematical programming" | 205 | 5 |
| 91 | Game theory, economics, social and behavioural sciences | | |
| | "game theory" | 55 | 2 |
| | "economics" | 508 | 5 |
| | "social science" | 87 | 2 |
| | "behavioural science" | 19 | 0 |
| 92 | Biology and other natural sciences | | |
| | "biology" | 501 | 14 |
| | "natural science" | 257 | 3 |
| 93 | Systems theory; control (including optimal control) | | |
| | "systems theory" or "theory of systems" | 140 | 0 |
| | "control theory" | 302 | 5 |
| | "optimal control" | 812 | 48 |
| 94 | Information and communication, circuits | | |
| | "information and communication" | 49 | 0 |
| | "circuit" | 746 | 44 |
| 97 | Mathematics education | | |

| | | |
|---|---|---|
| mathematics education or "math education", or "maths education" or "mathematical education." | 244 | 0 |

## 4. Discussion

The bibliometric analysis of published literature is essential to understand a research field from a bird's eye view: to understand what topics were popular and unpopular at the time, the significance of the works and authors, the links between them and other fields, and to find new directions for future research. Unfortunately, not all the manuscripts humanity publishes are always available in one academic database such as Scopus. Still, they are distributed over the internet in many databases, so extracting data can take a long time. Perhaps one day in the future, we will have a single tool that will allow academics to easily access all the knowledge that humanity has generated throughout its existence, summarise it quickly, translate it between languages, and use integrated AI agents as intellectual assistants.

The bibliometric research on the Bernstein polynomials allowed us to learn new visualisation tools and see one database's advantages over another. For example, Google Scholar indexes many more works than Scopus, but extracting bibliometric data is more complex and requires the use of software such as Publish and Perish[16]. I have already discussed this in the past, and I would like to take the opportunity to remind readers of my short article published in the leading English daily in South Korea [23].

## 5. Conclusions

In this work, we carried out a bibliometric analysis of Scopus data of documents containing the keyword "Bernstein polynomial" in all fields. The data extracted from Scopus allowed us to understand in which fields authors work, the dynamics and growth of publications, the most popular source titles, the dynamics of the number of publications in the four most popular journals, the percentage of documents by country, the leading universities whose authors mention Bernstein polynomials in their works, the most active academics, co-authorship between academics worldwide, funding agencies that often support research, the most cited works that study Bernstein polynomials. We also understood better that specific fields are strongly or weakly related to Bernstein polynomials, which can be seen in the visualisations. More importantly, we highlighted several areas with a weak connection to Bernstein polynomials and initiated a chat with AI chatbots, which generated some suggestions.

[16] https://harzing.com/resources/publish-or-perish